\documentclass[MSNbibl,number,citesort,seceqn,dvips]{arxbj}

\aid{0}
\volume{19}
\issue{3}
\pubyear{2013}
\firstpage{780}
\lastpage{802}
\doi{10.3150/11-BEJ416} %kopijuoti is PTS

\makeatletter
\newtheorem{theorem}{Theorem}[section]
\newremark{remark}[theorem]{Remark}
\newtheorem{lemma}[theorem]{Lemma}
\newtheorem{proposition}[theorem]{Proposition}
\newremark{example}[theorem]{Example}
\newtheorem{corollary}[theorem]{Corollary}

\def\join(#1;#2){[#1 \bullet\!\!\!- #2]}
\def\joinp(#1;#2){\Biggl[#1 \bullet\!\!\!- #2\Biggr]}
\def\joind(#1;#2){\bigl[#1 \bullet\!\!\!- #2\bigr]}
\def\N{\mathbb{N}}
\renewcommand{\mid}{\vert}
\makeatother

\begin{document}
\begin{frontmatter}

\title{Unimodularity for multi-type Galton--Watson~trees}
\runtitle{Multi-type Galton--Watson trees}

\begin{aug}
%%%% inicialai - be tarpu
\author{\fnms{Serdar} \snm{Altok}\corref{}\ead[label=e1]{serdar.altok@boun.edu.tr}}% \and
\runauthor{S. Altok} %% auto
\address{Department of Mathematics, Bogazici University, Bebek,
Istanbul 34342, Turkey.\\\printead{e1}}
\end{aug}

% HISTORY:
\received{\smonth{3} \syear{2010}}
\revised{\smonth{3} \syear{2011}}

% ABSTRACT
%
\begin{abstract}
Fix $ n \in\N$. Let $\mathbf{T}_n$ be the set of rooted trees
$(T,o)$ whose vertices are labeled by elements of $\{1,\dots,n\}$.
Let $\nu$ be a strongly connected multi-type Galton--Watson measure.
We give necessary and sufficient conditions for the existence of a
measure $\mu$ that is reversible for simple random walk on $\mathbf{T}_n$ and has the property that given the labels of the root and its
neighbors, the descendant subtrees rooted at the neighbors of the
root are independent multi-type Galton--Watson trees with conditional
offspring distributions that are the same as the conditional
offspring distributions of $\nu$ when the types are $\nu$ are
ordered pairs of elements of $[n]$. If the types of $\nu$ are given
by the labels of vertices, then we give an explicit description of
such $\mu$.
\end{abstract}

% KEYWORDS

\end{frontmatter}

%s1 ###
%s1 #&#
\section{Introduction}

Consider simple random walk on a multi-type Galton--Watson tree $T$
with distribution $\nu$. This induces a simple random walk on rooted
trees where the root represents the location of the original random
walker walking on the unrooted tree $T$. A stationary measure $\mu$
for this random walk can be useful in many ways such as calculating
the speed of simple random walk (see~\cite{augmented,Takacs}). If in addition $\mu$
is reversible for simple random walk, one has more results through
the connection between reversibility of simple random walk and
unimodularity (see~\cite{aldous}), which we will explain below.

The case when $\nu$ is a multi-type Galton--Watson measure with
deterministic conditional offspring distributions is well
understood: In~\cite{altok}, we showed that reversible measures for simple
random walk on rooted trees whose descendant subtrees are
deterministic multi-type Galton--Watson trees exist only for
Galton--Watson trees that are descendant subtrees of a cover of a
finite, connected, undirected graph. Such Galton--Watson trees are
periodic trees, that is, covers of finite, strongly connected, directed
graphs. Stationary measures for simple and biased random walk on
trees with periodic subtrees were previously studied in~\cite{Takacs}.

We first define covers of undirected and directed graphs. Let $H$ be
a finite connected undirected graph. Label the vertices of $H$ such
that any two vertices $x,y$ have the same label if and only if the
two rooted graphs $H$ rooted at $x$ and $H$ rooted at $y$ are rooted
isomorphic. Fix a vertex $x$ of $H$ and consider the unrooted tree
$T_x$ whose vertices are finite paths on $H$ that start from $x$ and
do not backtrack, where two vertices are connected by an edge if one
is the extension of the other by one edge. It is easy to see that
for all $x,y \in V(H)$, we have $T_x$ is isomorphic to $T_y$ (when
both are unrooted). Any tree that is isomorphic to $T_x$ for some $x
\in V(H)$ is called a \textit{cover} of $H$. The labeling on $H$ lifts
to a labeling of the vertices of the covers of $H$ by labeling each
vertex of a cover of $H$ by the label of the last vertex of the path
on $H$ that represents that vertex.

Similarly one can define covers of a finite, connected, directed
graph $G$. Label the vertices of $G$ such that any two vertices
$x,y$ have the same label if and only if the two rooted graphs $G$
rooted at $x$ and $G$ rooted at $y$ are rooted isomorphic. For each
vertex $x$ of $G$, we get a cover $T_x$ whose vertices are finite
directed paths in $G$ that start at $x$ and two vertices of $T_x$
are connected by an (undirected) edge if the path corresponding to
one vertex is the extension of the path corresponding to the other
vertex by a directed edge. The covers of $G$ are not necessarily
isomorphic.

Now, we explore a connection between covers of undirected and
directed graphs, namely, we show that the descendant subtrees of
covers of an undirected graph are covers of a directed graph. Let
$H$ and $T_x$ for $x \in V(H)$ be as above. Since the vertices of
$T_x$ are finite paths in $H$ starting at $x \in V(H)$, there is a
vertex of $T_x$ that corresponds to the path consisting of just the
vertex $x \in V(H)$. Call this vertex of $T_x$, $x$, as well.
Consider $T_x$ as rooted at $x \in V(T_x)$. Delete $x \in V(T_x)$,
hence the edges incident to it, in $T_x$. The remaining trees, the
descendant subtrees of the neighbors of $x \in V(T_x)$, are covers
of the same directed graph. To see this remember that the vertices
of $T_x$ are labeled by the label of the last vertex of the
corresponding path on $H$. Thus, if $y \in V(T_x)$ has label $i$ and
the parent of $y$, the neighbor of $y$ on the path to $x$, has label
$j$, it means that the last edge of the path in $H$ corresponding to
$y$ has endpoints labeled $j$ and $i$ with $i$ being the label of
the last vertex of this path. Since the paths do not backtrack any
vertex of $T_x$ that has label $i$ and whose parent has label $j$,
has the same descendant tree. Now, we can describe the directed
graph $G$, whose covers are exactly the descendant subtrees of
$T_x$. Vertices of $G$ are ordered pairs of labels of $H$ and for
any labels $i,j,k$, $(j,i)$ is connected to $(i,k)$ by $m$ directed
edges if a vertex of $H$ labeled $i$ has a neighbor labeled $j$ and
$m$ other neighbors that are labeled~$k$.

Let $(T,o)$ denote the tree $T$ rooted at $o \in V(T)$. Then for any
rooted cover $(T,o)$ of a finite connected undirected graph $H$,
when $o$ is deleted, the remaining trees are deterministic
multi-type Galton--Watson trees, as explained in the previous
paragraph, where the type of a vertex is determined by a one-to-one
function of the label of the vertex and of its parent. The uniform
measure on the vertices of $H$ biased by the degree of the vertex
lifts to a measure $\mu$ on rooted covers of $H$ that is reversible
for simple random walk. This also establishes the existence of a
unimodular measure with the same support since reversibility of
$\mu$ is equivalent to the unimodularity of the probability measure
obtained by biasing $\mu$ by the reciprocal of the degree of the
root. We will obtain our results in the setting of reversible
measures on rooted trees.

This relationship for (single-type, non-deterministic) Galton--Watson
trees was established in~\cite{augmented}. Let $\nu$ be a Galton--Watson measure
with offspring distribution $\{p_k\}_{k \in\N}$. Start with one
vertex, the root, that has $k+1$ neighbors with probability $p_k$,
each of which has independent Galton--Watson descendant\vadjust{\goodbreak} subtrees with
offspring distribution~$\nu$. This is the Augmented Galton--Watson
measure and is reversible for simple random walk on rooted trees.
Biasing this measure by the reciprocal of the degree of the root,
one gets the unimodular Galton--Watson measure.

Now, let $\nu$ be a multi-type Galton--Watson measure. Motivated by
the results on deterministic multi-type Galton--Watson trees and
non-deterministic single-type Galton--Watson trees, we try to find
reversible measures $\mu$ for simple random walk on labeled rooted
trees with the following property: Conditioned on the labels of the
root and its neighbors the subtrees rooted at the neighbors of the
root are independent (non-deterministic) multi-type Galton--Watson
subtrees whose conditional offspring distributions are the same as
the conditional offspring distributions of $\nu$.

The deterministic case shows that if $\mu$ satisfies the property
above, then the types of $\nu$ are not necessarily given by the
labels of vertices. Therefore, we introduce relabeling functions on
rooted labeled trees.

For $n \in\N$, let $\mathbf{T}_n$ be the set of rooted trees $(T,o)$
whose vertices are labeled by elements of $\{1,\dots,n\}$. Let $f$
be a function from $\mathbf{T}_n$ to rooted trees whose vertices are
labeled except for the root, with the following properties:
$V(f(T,o))=V(T,o),E(f(T,o))=E(T,o)$ and $f(T,o)$ is rooted at $o$.
Let $\mu$ be a measure on $\mathbf{T}_n$. Suppose that conditioned on
the labels of the root and its neighbors, after relabeling the
vertices by $f$ the descendant subtrees of the neighbors of the root
are independent Galton--Watson trees with conditional offspring
distributions the same as those of $\nu$. In this case, we write $\mu
\sim_f \nu$ and the label of a vertex after relabeling gives the
$\nu$-type of that vertex.

We consider two cases: In the first case, the new label of a vertex
is given by a one-to-one function of the label of the vertex and of
its parent. That is, if the labels of $x \in V(f(T_1,o_1))$ and
$y\in V(f(T_2,o_2))$ are the same, then $x\in V(T_1,o_1)$ and $y\in
V(T_2,o_2)$ have the same label and their parents have the same
label as well. For a given $\nu$, we give necessary and sufficient
conditions on the conditional offspring distributions of $\nu$ for
the existence of a measure $\mu$ such that $\mu\sim_f \nu$ where
$f$ is as above and $\mu$ is reversible for simple random walk (see
Theorem~\ref{t2} and Example~\ref{xex1}). We show that among all the
measures that
have a given conditional offspring distribution, there is at most
one measure $\mu$ such that $\mu\sim_f \nu$ for some strongly
connected Galton--Watson measure $\nu$ and $\mu$ is reversible and we
parametrize all such $\mu$ (see Proposition~\ref{pdimension} and
Example~\ref{xex2}).

We also study the case when there is no relabeling, that is, for all
$(T,o)$ and all $x\in V(T,o)$, the label of $x$ in $(T,o)$ is the
same as the label of $x$ in $f(T,o)$. In this case, given the degree
and the type of the root, the conditional offspring distributions of
the root with respect to $\mu$ and $\nu$ are multinomial and the
parameters for the conditional offspring distributions do not depend
on the degree of the root (see Theorem~\ref{tnorelabeling} and Example
\ref{xex3}).

%%%%%%%%%%%%%%%%%%%%%%%%%%%%%%%%%%%%%%%%%%%%%%%%%%%%
%%%%%%%%%%%%%%%%%%%%%%%%%%%%%%%%%%%%%%%%%%%%%%%%%%%%
%%% %%%%%%%%%%%%%%
%%% SECTION DEFINITIONS AND NOTATIONS %%%%%%%%%%%%%%
%%% %%%%%%%%%%%%%%
%%%%%%%%%%%%%%%%%%%%%%%%%%%%%%%%%%%%%%%%%%%%%%%%%%%%
%%%%%%%%%%%%%%%%%%%%%%%%%%%%%%%%%%%%%%%%%%%%%%%%%%%%

%s2 ###
%s2 #&#
\section{Definitions and notations}

A rooted tree $(T,o)$ is a connected graph $T$ with no cycles that
has a distinguished vertex~$o$. For any $x \in V(T,o)$ other than $o$,
we call the
neighbor of $x$ on the unique simple path between $x$ and $o$ the
\textit{parent} of $x$. The remaining neighbors of $x$ are called \textit{children} of $x$. The neighbors of the root are
called children of the root.

We shall work with rooted isomorphism classes of rooted labeled
trees: For all $n \in\N$ let $[n]:=\{1,\dots,n\}$. A labeled rooted
tree is a rooted tree whose vertices are labeled by elements of some
fixed $[n]$. Two rooted labeled trees $(T,o),(T',o')$ are called
\textit{rooted isomorphic} if there is a bijection between the
vertices of the two trees that maps $o$ to~$o'$, preserves adjacency
and any two vertices of $T$ that have the same label are mapped to
vertices that have same label and conversely. We shall use $(T,o)$
to denote the rooted isomorphism class of the labeled rooted tree
$(T,o)$ as well. For all $n$, let $\mathbf{T}_n$ denote the set of
rooted isomorphism classes of rooted trees labeled by elements of
$[n]$.

Fix $n \in\N$. Let $\nu$ be a Galton--Watson measure
whose types are elements of $[n]$. We say that
$\nu$ is \textit{strongly connected} if for all $i,j \in[n]$, the root
is of type $i$ with positive probability and a tree has a vertex of
type $j$ with positive probability given that the root has type $i$.

Let $\mathbf{S}_n:=\{(c_1,\dots,c_n)\dvt c_i \in\N\cup{0} \}$. For $\mathbf{c}=(c_1,\dots,c_n) \in\mathbf{S}_n$,
let $\mathbf{c}_j:=(c_1,\dots,c_j-1,\allowbreak\dots,c_n)$, $\mathbf{c}^k:=(c_1,\dots,c_k+1,\dots,c_n)$, $\mathbf{c}_j^k:=(\mathbf{c}_j)^k$ and
$|\mathbf{c}|:=\sum_i c_i$. We use numbers in parenthesis to index a
set of vectors. If $\mathbf{c}(1),\dots,\mathbf{c}(m)$ are $m$ vectors,
then $c(i)_p$ denotes the $p$th coordinate of the vector $\mathbf{c}(i)$.

We say that a vertex has $\mathbf{c}$ neighbors if it has exactly $c_j$
neighbors labeled $j$ for every $j \in[n]$. For all $i \in\N$, let
$N_i$ be the set of trees rooted at a vertex labeled $i$. Let
$N_i(\mathbf{c})$ be the set of trees whose roots are labeled $i$ and
have $\mathbf{c}$ neighbors.

For a probability measure $\sigma$ on $\mathbf{T}_n$ and $i \in[n]$,
let $\sigma(i):=\sigma(N_i)$.

For any measure $\sigma$ on $\mathbf{S}_n$, we say that $\sigma$ has
the \textit{multinomial} distribution with \textit{parameters}
$p_1,\dots,p_n \in[0,1]$ satisfying $\sum_{i=1}^n p_i=1$ if
for some $l \in\N$, we have $\sigma(\mathbf{c})=\frac{|\mathbf{c}|!}{\prod_{k\in[n]} c_k!}\prod_{k
\in[n]} p_k^{c_k}$ for all $\mathbf{c} \in\mathbf{S}_n$ with $|\mathbf{c}|=l$.

For all $k \in\N,\mathbf{c}, \mathbf{d} \in\mathbf{S}_n$ and
$i_1,i_2,\ldots,i_k \in[n]$, define $B_{\mathbf{c},i_1,i_2,\dots,i_k,
\mathbf{d}}$ to be the set of trees whose roots are labeled $i_1$ and whose
roots have $\mathbf{c}^{i_2}$ children and at least one child of the
root that is labeled $i_2$ has at least one child labeled $i_3$,
which has at least one child labeled $i_4$ and so on up to a vertex
labeled $i_k$, which has $\mathbf{d}$ children. We also define
\begin{eqnarray*}
B_{\mathbf{c},i_1,i_2,\dots,i_k}&:=&\bigcup_\mathbf{d}B_{\mathbf{c},i_1,i_2,\dots,i_k,\mathbf{d}},
\\
B_{i_1,i_2,\dots,i_k,\mathbf{d}}&:=&\bigcup_\mathbf{c}B_{\mathbf{c},i_1,i_2,\dots,i_k,\mathbf{d}},
\\
B_{i_1,i_2,\dots,i_k}&:=&\bigcup_{\mathbf{c},\mathbf{d}}B_{\mathbf{c},i_1,i_2,\dots,i_k,\mathbf{d}} .
\end{eqnarray*}

For a probability measure $\rho$ on $\mathbf{S}_n$, we define the
following sets:
\begin{eqnarray*}
A(\rho)&:=&\bigl\{j \in[n]\dvt \exists\mathbf{c} \mbox{ such that }
c_j>0,\rho (\mathbf{c})>0\bigr\},
\\
D(\rho)&:=&\bigl\{d\dvt \exists\mathbf{c} \mbox{ with } |\mathbf{c}|=d,\rho(
\mathbf{c})>0 \bigr\}.
\end{eqnarray*}

For a probability measure $\sigma$ on $\mathbf{T}_n$ and $i \in[n]$,
let $\sigma_i$ denote the conditioning of $\sigma$ on $N_i$. Let
$\sigma_i(\mathbf{c}):=\sigma_i(N_i(\mathbf{c}))$. This way we interpret
$\sigma_i$ as a probability measure on $\mathbf{S}_n$.

Thus, the neighbors of a vertex $x$ of $(T,o) \in\mathbf{T}_n$ whose
label is $i$ are labeled by elements of $A({\mu_i})$ $\mu$-a.s. and
the number of neighbors of $x$ is in $D(\mu_i)$ $\mu$-a.s. We write
$\mathbf{c} \in A(\rho)$ if $c_j=0$ for every $j\notin A(\rho)$. Hence,
$\rho(\mathbf{c})>0$ implies $\mathbf{c} \in A(\rho)$.

%%%%%%%%%%%%%%%%%%%%%%%%%%%%%%%%%%%%%%%%%%%%%%%%%%%%%%%%%%%%%%%%%
%%%%%%%%%%%%%%% MAIN RESULTS %%%%%%%%%%%%%%%%%%%%
%%%%%%%%%%%%%%%%%%%%%%%%%%%%%%%%%%%%%%%%%%%%%%%%%%%%%%%%%%%%%%%%%
%%%%%%%%%%%%%%%%%%%%%%%%%%%%%%%%%%%%%%%%%%%%%%%%%%%%%%%%%%%%%%%%%

%s3 ###
%s3 #&#
\section{Main results}

Remember that the relabeling function $f$ is a function from rooted
labeled trees to rooted trees whose vertices are labeled except for
the root such that the vertices, edges and the root of $f(T,o)$ is
the same as those of $(T,o)$. So if we forget about the labels of
$(T,o)$ and $f(T,o)$, they are the same rooted tree. Fix $f$ such
that for all $(T,o),(T',o')\in\mathbf{T}_n$ and $x \in V(T,o),y\in
V(T',o')$ the labels of $x$ in $f(T,o)$ and $y$ in $f(T',o')$ are
the same if and only if the labels of $x$ in $(T,o)$ and $y$ in
$(T',o')$ are the same and the parent of $x$ in $(T,o)$ has the same
label as the parent of $y$ in $(T',o')$. Thus, if $x \in V(T,o)$ has
label $j$ and its parent has label $i$, then we denote the label of
$x$ in $f(T,o)$ by $(i,j)$. We write $\nu_{i,j}$ for $\nu_{(i,j)}$
and $\nu_{i,j}(\mathbf{c})$ for
$\nu_{i,j}((j,c_1),\dots,(j,c_n))=\nu_{i,j}(N_{(i,j)}((j,c_1),\dots
,(j,c_n)))$.

For this section, $\mu$ denotes
a measure on $\mathbf{T}_n$ such that $\mu(i)>0$ for all $i\in[n]$ and
$\nu$ denotes a Galton--Watson measure whose
types are ordered pairs of elements of $[n]$. The relabeling function
$f$ is fixed as above and
we write $\mu\sim\nu$ instead of $\mu\sim_f \nu$.

Our first goal is to find necessary and sufficient conditions on $\nu$
for the existence of a reversible measure $\mu$ for simple random walk
such that $\mu\sim\nu$.
We start by studying properties of $\mu$ and $\nu$ when $\mu\sim\nu$.

%re3.1 #&#
\begin{remark}\label{r1}
Assume $\mu\sim\nu$. For all $i_1,\dots,i_k \in[n]$
and for all $\mathbf{c} \in\mathbf{S}_n$, if $\mu(B_{i_1,\dots,i_k,\mathbf{c}})>0$, then the following two properties hold:
\begin{enumerate}[(ii)]

\item[(i)] for all $j \in[2,n-1],(i_j,i_{j+1}) \in A(\nu_{i_{j-1},i_j})$,

\item[(ii)] $\nu_{i_{k-1},i_k}(\mathbf{c})>0$.
\end{enumerate}
\end{remark}

%%%%%%%%%%%%%%%
%%% LEMMA 1 %%%
%%%%%%%%%%%%%%%

%le3.2 #&#
\begin{lemma}\label{l1}
Assume that $\mu\sim\nu$ and that $\mu$ is stationary
for simple random walk. Then for all $\mathbf{c} \in\mathbf{S}_n$ and
$i,j \in[n]$ such that $c_j>0$,
%
%e3.1 #&#
\begin{equation}
\label{emunu} \mu\bigl(N_i(\mathbf{c})\bigr)>0 \mbox{ (or
equivalently $\mu_i(\mathbf{c})>0$)}\quad \iff\quad
\nu_{j,i}(\mathbf{c}_j)>0.
\end{equation}
In this case, we have
%
%e3.2 #&#
\begin{equation}
\label{eijji} \forall i,j,k \in[n] \mbox{ and } \mathbf{c} \in
\mathbf{S}_n \mbox{ with } c_j,c_k>0 \qquad
\nu_{k,i}(\mathbf{c}_k)>0 \quad\iff\quad\nu_{j,i}(
\mathbf{c}_j)>0.
\end{equation}
\end{lemma}

\begin{pf}
Fix $i,j$ and $\mathbf{c}$ with $c_j>0$.

Assume $\mu(N_i(\mathbf{c}))>0$. Then $\mu(B_{\mathbf{c}_j,i,j})>0$. A
random walker at $(T,o)\in B_{\mathbf{c}_j,i,j}$ can walk to a tree in
$B_{j,i,\mathbf{c}_j}$ with positive probability. Since $\mu$ is
stationary and $\mu(B_{\mathbf{c}_j,i,j})>0$, we have $\mu(B_{j,i,\mathbf{c}_j})>0$. By Remark~\ref{r1} we have $\nu_{j,i}(\mathbf{c}_j)>0$.

Conversely, assume $\nu_{j,i}(\mathbf{c}_j)>0$. Since $\nu$ is strongly
connected, there exist $i_1,\dots,i_{k-1}=j,i_k=i$ such that
$\mu(B_{i_1,i_2,\dots,j,i,\mathbf{c}_j})>0$. A random walker at
$(T,o)\in B_{i_1,i_2,\dots,j,i,\mathbf{c}_j}$ can walk to a tree in
$B_{\mathbf{c}_j,i,j,...,i_2,i_1}$ with positive probability. Since
$\mu$ is stationary, this gives $\mu(B_{\mathbf{c}_j,i,j,\dots,i_2,i_1})>0$, which implies $\mu(N_i(\mathbf{c}))=\mu(B_{\mathbf{c}_j,i,j})>0$.

The statement in (\ref{eijji}) is a direct corollary of the first
statement of the lemma.\vadjust{\goodbreak}
\end{pf}

Let $\mu$ be a measure on $\mathbf{T}_n$ and $\nu$ be a multi-type
Galton--Watson measure with $\mu\sim\nu$.
For $A \in\mathbf{T}_n$, let $\mathbf{P}((T,o) \to A)$ denote the
probability that a simple random walker that starts at $(T,o)$ moves
to some $(T,o') \in A$ in one step. Then $\mu$ is reversible if and
only if for all measurable sets $A,B \subset\mathbf{T}_n$
%
%e3.3 #&#
\begin{equation}
\label{ereversible} \int_A\mathbf{P}\bigl((T,o) \to B
\bigr) \,\mathrm{d}\mu(T,o)=\int_B\mathbf{P}\bigl((T,o) \to
A\bigr)\, \mathrm{d}\mu(T,o).
\end{equation}

We first prove a lemma that says $\mu$ is reversible if
and only if the reversibility equations (\ref{ereversible})
hold for a special family of
measurable sets.

%%%%%%%%%%%%%%%%%%%%%%%%
%%% LEMMA REVERSIBLE %%%
%%%%%%%%%%%%%%%%%%%%%%%%
%
%le3.3 #&#
\begin{lemma}\label{lreversible}
Let $\mu$ be a measure on $\mathbf{T}_n$. Let $\nu$
be a multi-type Galton--Watson measure. Assume that $\mu\sim\nu$
and (\ref{emunu}) holds. Then the following statements are
equivalent:
\begin{enumerate}[(iii)]
\item[(i)] $\mu$ is reversible,

\item[(ii)] for all $i,j \in[n]$ and $\mathbf{c},\mathbf{d}\in\mathbf{S}_n$ with
$c_j,d_i>0$ we have
%
%e3.4 #&#
\begin{equation}
\label{e2} \mu(i)\mu_i(\mathbf{c})\nu_{i,j}(
\mathbf{d}_i)\frac{c_j }{|\mathbf{c}|} =\mu(j)\mu_j(\mathbf{d})
\nu_{j,i}(\mathbf{c}_j)\frac{d_i }{|\mathbf{d}|},
\end{equation}

\item[(iii)] for all $i,j \in[n]$ and $\mathbf{c},\mathbf{d}\in\mathbf{S}_n$
%
%e3.5 #&#
\begin{equation}
\label{elocalreversible} \int_{N_i(\mathbf{c})}{\mathbf{P}\bigl((T,o) \to
N_j(\mathbf{d})\bigr)\, \mathrm{d}\mu(T,o)} =\int
_{N_j(\mathbf{d})}{\mathbf{P}\bigl((T,o) \to N_i(\mathbf{c})
\bigr)\, \mathrm{d}\mu(T,o)}.
\end{equation}
\end{enumerate}
\end{lemma}

%%%%%%%%%%%%%%%%%%%%%%%%%%%%%%%%%
%%% PROOF OF LEMMA REVERSIBLE %%%
%%%%%%%%%%%%%%%%%%%%%%%%%%%%%%%%%

\begin{pf}
Statement (i) is equivalent to (\ref{ereversible}), so our
task is to show that (\ref{ereversible}), (\ref{e2}) and (\ref
{elocalreversible}) are equivalent.

%%%%%%%%%%%%%%%%%%%%%%%%%%
%%% DETERMINISTIC PART %%%
%%%%%%%%%%%%%%%%%%%%%%%%%%

If $\nu$ is a deterministic multi-type Galton--Watson measure, then
any measurable set is a finite disjoint union of the sets $N_i(\mathbf{c})=N_i$ up to sets of measure zero. Hence, (\ref{ereversible}) and
(\ref{elocalreversible}) are equivalent. On the other hand when $\mu$
is deterministic, (\ref{elocalreversible}) is equivalent to
\[
\mu\bigl(N_i(\mathbf{c})\bigr)\frac{c_j }{|\mathbf{c}|}=\mu
\bigl(N_j(\mathbf{d})\bigr)\frac{d_i }{
|\mathbf{d}|}
\]
since $\mathbf{P}((T,o) \to N_j(\mathbf{d}))=\frac{c_j }{|\mathbf{c}|}$ for all
(actually the unique) $(T,o) \in N_i(\mathbf{c})$. Both sides of
(\ref{e2}) and (\ref{elocalreversible}) are zero unless the root has
label $i$ and has neighbors $\mathbf{c}$ with $c_j>0$ with positive
probability and the root has
label $j$ and has neighbors $\mathbf{d}$ with positive probability. Since
$\mu(i)\mu_i(\mathbf{c})=\mu(N_i(\mathbf{c}))$, $\mu(j)\mu_j(\mathbf{d})=\mu(N_j(\mathbf{d}))$ and $\nu_{i,j}
(\mathbf{d}_i)=\nu_{j,i}(\mathbf{c}_j)=1$, we conclude that (\ref{elocalreversible}) is equivalent to
(\ref{e2}).

%%%%%%%%%%%%%%%%%%%
%%% RANDOM PART %%%
%%%%%%%%%%%%%%%%%%%

Let $\nu$ be a non-deterministic multi-type Galton--Watson measure.
We first show (\ref{ereversible}) is equivalent to condition (ii). We
emulate the proof of Theorem 3.1 in~\cite{augmented}.

Given two rooted trees $(T,o),(T',o')\in\mathbf{T}_n$, define $[(T,o)
\bullet\!\!\!-(T',o')]$ to be the tree formed by joining disjoint
copies of $(T,o)$ and $(T',o')$ by an edge between $o$ and $o'$ and
rooting the new tree at $o$. For sets $F,F' \subset\mathbf{T}_n,$
write
\[
\join(F;F'):=\bigl\{\bigl[(T,o) \bullet\!\!\!- {
\bigl(T',o'\bigr)}\bigr]\dvt (T,o)\in F,
\bigl(T',o'\bigr)\in F'\bigr\}.
\]

For trees $(T_1,o_1),\dots,(T_m,o_m)\in\mathbf{T}_n$ and $i\in[n]$, let
$\join(i;\bigvee_{s=1}^m {(T_s,o_s)})$ denote the tree formed by
joining the vertices $o_s$ by single edges to a new vertex labeled
$i$, the latter being the root of the new tree. Furthermore, for
sets $F_1,\dots,F_m$, write
\[
\joinp(i;\bigvee_{s=1}^m
F_s):=\Biggl\{\joinp(i;\bigvee
_{s=1}^m {(T_s,o_s)})
\dvt (T_s,o_s)\in F_s\Biggr\}.
\]

For any rooted tree $(T,o)$ and $x \in V(T)$, the \textit{distance} of
$x$ to $o$ is the number of edges on the simple path between $x$ and
$o$. For any rooted tree $(T,o)$ and for a positive integer~$h$, the
restriction of $(T,o)$ to the first $h$ levels is the finite tree
induced by the vertices whose distance to $o$ is less than or equal
to $h$. The \textit{height} of a finite rooted tree $(T,o)$ is the
maximum of the distances of the vertices of $T$ to $o$.

For a finite tree $(T,o)$ of height $h\in\N$, let $B_{(T,o)}$ be
the set of rooted trees whose restriction to the first $h$ levels is
isomorphic to $(T,o)$. The sets $B_{(T,o)}$ generate the
$\sigma$-field. Let $K$ be the set of all $B_{(T,o)}$. Then any
element of $K$ is of the form
$\join(B_{(T,o)};B_{(T',o')})=B_{\join((T,o);{(T',o')})}$ where
%
%e3.6 #&#
\begin{equation}
\label{etctd} (T,o):=\joinp(i;\bigvee_{s=1}^l
{(T_s,o_s)})\quad\mbox{and}\quad
\bigl(T',o'\bigr):=\joinp(j;\bigvee
_{s=1}^m {\bigl(T'_s,o'_s
\bigr)})
\end{equation}
for finite trees $(T_1,o_1),\dots,(T_l,o_l)$ and
$(T'_1,o'_1),\dots,(T'_m,o'_m)$.

We now show that (\ref{ereversible}) holds for all sets of the form
$A=B_{\join((T,o);{(T',o')})} \in K$ and
$B=B_{\join((T',o');{(T,o)})}$ if and only if (\ref{e2}) holds. Let
$l(x)$ denote the label of a vertex $x$. Fix
$\join((T, o);{(T',o')})$. Let $h$ be the height of
$\join((T,o);{(T',o')})$. Without loss of generality, assume $(T,o)$
has height $h$. Let $\mathbf{c},\mathbf{d} \in\mathbf{S}_n$ be such that for
all $k$ we have
\[
c_k=\bigl|\bigl\{x \in V\bigl(\joind((T,o);{\bigl(T',o'
\bigr)})\bigr)\dvt x \sim o,l(x)=k\bigr\}\bigr|
\]
and
\[
d_k=\bigl|\bigl\{x \in V\bigl(\joind(\bigl(T',o'
\bigr);{(T,o)})\bigr)\dvt x \sim o',l(x)=k\bigr\}\bigr|.
\]
Let $j=l(o')$ and $i=l(o)$. Then we have $c_j>0$ and $d_i>0$. Define
$c^*_j$ by
\[
c^*_j:=\bigl|\bigl\{(T_s,o_s)\dvt
(T_s,o_s) \cong\bigl(T',o'
\bigr),s \in[l]\bigr\}\bigr|+1.
\]
Then the probability that the chain moves from an element
of $A$ to an element of $B$ is~$\frac{c^*_j }{|\mathbf{c}|}$.

Let $(T^*,o)$ be the restriction of $(T,o)$ to the first $h-2$
levels. Let $d^*_i$ be 1 plus the number of $o'_k \sim o' \in V(T')$
for $k=1,\dots,m$ such that $(T^*,o)$ is isomorphic to
$(T'_k,o'_k)$.

Let $D,D'$ be such that
%
%e3.7 #&#
\begin{equation}
\label{eji} \nu_{j,i}(B_{(T,o)})=\nu_{j,i}(
\mathbf{c}_j){c_j-1 \choose c^*_j-1}
\nu_{i,j}(B_{(T',o')})^{c^*_j-1}D
\end{equation}
and
%
%e3.8 #&#
\begin{equation}
\label{eij} \nu_{i,j}(B_{(T',o')})=\nu_{i,j}(
\mathbf{d}_i){d_i-1 \choose d^*_i-1}
\nu_{j,i}(B_{(T^*,o)})^{d^*_i-1}D'.
\end{equation}

We first calculate $\mu(B_{\join((T,o);{(T',o')})})$. In (\ref{eji}),
the term $D$ takes care of all the descendant subtrees of the
neighbors of $o$ that are not labeled $j$ and also the descendant
subtrees of the $c_j-1 - (c^*_j-1)=c_j - c^*_j$ neighbors of $o$
that are labeled $j$ and have descendant subtrees whose restriction
to the first $h-1$ levels are not $(T',o')$. To calculate
$\mu(B_{\join((T,o);{(T',o')})})$, we choose $c_j-c^*_j$ of the $c_j$
neighbors of the root that are labeled $j$ to have descendant
subtrees whose restrictions to the first $h-1$ levels are not
$(T',o')$. We have
\[
\mu(B_{\join((T,o);{(T',o')})})=\mu(i)\mu_i(\mathbf{c}){c_j
\choose c_j-c^*_j}\nu_{i,j}(B_{(T',o')})^{c^*_j}D,
\]
which by (\ref{eij}) is equal to
\[
\mu(i)\mu_i(\mathbf{c}){c_j \choose c^*_j}
\nu_{i,j}(B_{T',o'})^{c^*_j-1}D\nu_{i,j}(
\mathbf{d}_i){d_i-1 \choose d^*_i-1}
\nu_{j,i}(B_{(T^*,o)})^{d^*_i-1}D'.
\]

Then the left-hand side of (\ref{ereversible}) for
$A=B_{\join((T,o);{(T',o')})}$ and $B=B_{\join((T',o');{(T,o)})}$ is
equal to $\mu(B_{\join((T,o);{(T',o')})})\frac{c^*_j }{|\mathbf{c}|},$
which in turn is equal to
%
%e3.9 #&#
\begin{equation}
\label{ereversibleA} \mu(i)\mu_i(\mathbf{c})\nu_{i,j}(
\mathbf{d}_i)\frac{c_j }{|\mathbf{c}|}{c_j-1 \choose
c^*_j-1}\nu_{i,j}(B_{(T',o')})^{c^*_j-1}{d_i-1
\choose d^*_i-1} \nu_{j,i}(B_{(T^*,o)})^{d^*_i-1}DD'.
\end{equation}

Similarly, we calculate the right-hand side of (\ref{ereversible}).
For $k \in[d^*_i]$, let $B^k_{(T',o')}$ be the set of all trees
$(T'',o')$ such that
\begin{longlist}[(ii)]

\item [(i)] the restriction of $(T'',o')$ to the first $h-2$ levels is isomorphic
to $(T',o')$,

\item [(ii)] when $o' \in V(T'',o')$ is deleted, there are exactly $k-1$ trees
(out of the $|\mathbf{d}_i|=m$) whose restrictions
to the first $h$ levels are isomorphic to $(T,o)$.
\end{longlist}

So for any $(T'',o') \in B^k_{(T',o')}$, we have that $o'$ has
$d^*_i-1$ neighbors whose descendant trees' restrictions
to the first $h-2$ levels are
isomorphic to $(T^*,o)$ and exactly $k-1$ of them have descendant trees
whose restrictions to the first $h$ levels are isomorphic to $(T,o)$.
Let $E$ be defined by $\nu_{j,i}(B_{(T,o)})=\nu_{j,i}(B_{(T^*,o)})E$.
Now we have
\[
\mu\bigl(\join(B^k_{(T',o')};B_{(T,o)})\bigr)=
\mu(j)\mu_j(\mathbf{d}) {d_i \choose d^*_i}
\nu_{j,i}(B_{(T^*,o)})^{d^*_i}{d^*_i \choose
k}E^k(1-E)^{d^*_i-k}D'.
\]
Then the right-hand side of (\ref{ereversible}) for
$A=B_{\join((T,o);{(T',o')})}$ and
$B=\join(B^k_{(T',o')};B_{(T,o)})$ is equal to
$\mu(\join(B^k_{(T',o')};B_{(T,o)}))\frac{k }{|\mathbf{d}|}$. Since
$B_{\join((T',o');{(T,o)})}=\bigcup_{k=1}^{d^*_i}
\join(B^k_{(T',o')};B_{(T,o)})$ and the sets
\mbox{$\join(B^k_{(T',o')};B_{(T,o)})$} are disjoint, the right-hand side
of (\ref{ereversible}) for $A=B_{\join((T,o);{(T',o')})}$ and
$B=B_{\join((T',o');{(T,o)})}$ is given by
\begin{eqnarray*}
&&\sum_{k=1}^{d^*_i}\mu(j)\mu_j(
\mathbf{d}){d_i \choose d^*_i}\nu_{j,i}(B_{T^*,o})^{d^*_i}
{d^*_i \choose k}E^k(1-E)^{d^*_i-k}D'
\frac{k }{|\mathbf{d}|}
\\
&&\quad=\mu(j)\mu_j(\mathbf{d}){d_i \choose
d^*_i}\nu_{j,i}(B_{(T^*,o)})^{d^*_i}
d^*_iD'\frac{1 }{|\mathbf{d}|}E\sum
_{k=1}^{d^*_i}{d^*_i-1 \choose
k-1}E^{k-1}(1-E)^{d^*_i-1-(k-1)} .
\end{eqnarray*}
The sum on the right is equal
to 1. Since $\nu_{j,i}(B_{(T^*,o)})E=\nu_{j,i}(B_{T,o})$, the above
expression reduces to
\[
\mu(j)\mu_j(\mathbf{d})\frac{d_i }{|\mathbf{d}|}{d_i-1 \choose
d^*_i-1}\nu_{j,i}(B_{(T^*,o)})^{d^*_i-1}D'
\nu_{j,i}(B_{(T,o)}).
\]
Using (\ref{eji}) this is equal to
%
%e3.10 #&#
\begin{equation}
\label{ereversibleB} \mu(j)\mu_j(\mathbf{d})\frac{d_i }{|\mathbf{d}|}{d_i-1
\choose d^*_i-1}\nu_{j,i}(B_{(T^*,o)})^{d^*_i-1}D'
\nu_{j,i}(\mathbf{c}_j){c_j-1 \choose
c^*_j-1}\nu_{i,j}(B_{T',o'})^{c^*_j-1}D.
\end{equation}
Setting (\ref{ereversibleA}) equal to (\ref{ereversibleB}) we see
that when (\ref{emunu}) holds and
$A=B_{\join((T,o);(T',o'))}$ and $B=B_{\join((T',o');(T,o))}$,
(\ref{ereversible}) is equivalent to (\ref{e2}).

Since $K$ is a $\pi$-system, by the $\pi-\lambda$ theorem (\ref
{ereversible}) holds for all $\join(A;B)$ if and only if it holds for all
$\join(B_{(T,o)};B_{(T',o')})\in K$. Thus, our first claim is proved.

Next, we show that (\ref{elocalreversible}) is equivalent to (\ref
{e2}). If $(T,o) \in N_i(\mathbf{c})$, then a random walker at $(T,o)$
moves to a tree $(T',o') \in N_j(\mathbf{d})$ with probability $\frac{k}{
|\mathbf{c}|}$, where $k$ is the number of children of $o$ of type $j$
that have $\mathbf{d}_i$ children. Hence, we have
\begin{eqnarray*}
\mbox{LHS of }(\ref{elocalreversible}) &=&\sum_{k=1}^{c_j}
\mu\bigl(N_i(\mathbf{c})\bigr) \bigl(\nu_{i,j}(
\mathbf{d}_i)\bigr)^k\bigl(1-\nu_{i,j}(
\mathbf{d}_i)\bigr)^{c_j-k}{c_j \choose k}
\frac{k }{
|\mathbf{c}|}
\\
&=&\mu\bigl(N_i(\mathbf{c})\bigr)\nu_{i,j}(
\mathbf{d}_i)\frac{c_j }{|\mathbf{c}|}\sum_{k=1}^{c_j}
\bigl({\nu_{i,j}(\mathbf{d}_i)\bigr)^{k-1}\bigl(1-
\nu_{i,j} (\mathbf{d}_i)\bigr)^{c_j-k}{c_j-1
\choose k-1}}
\\
&=&\mu\bigl(N_i(\mathbf{c})\bigr)\nu_{i,j}(
\mathbf{d}_i)\frac{c_j }{|\mathbf{c}|}=\mu(i)\mu_i(\mathbf{c})
\nu_{i,j}(\mathbf{d}_i)\frac{c_j }{|\mathbf{c}|}.
\end{eqnarray*}
Making the same calculation for the right-hand side of (\ref
{elocalreversible}), we have that (\ref{elocalreversible}) is
equivalent to~(\ref{e2}).
\end{pf}

%%% END PROOF LEMMA REVERSIBLE %%%

Now assume that $\mu$ is reversible for simple random walk with $\mu
\sim\nu$.
Let $\mathbf{c},\mathbf{d},\mathbf{e} \in\mathbf{S}_n$ be such that
$c_j,d_i,e_i >0$ and $\mu_i(\mathbf{c}),\mu_j(\mathbf{d}),\mu_j(\mathbf{e})>0$.
By Lemma~\ref{l1}, we have that (\ref{eijji}) holds.
Then (\ref{elocalreversible}) holds for $A=N_i(\mathbf{c})$ and
$B=N_j(\mathbf{d})$
which reduces to (\ref{e2}). Similarly (\ref{elocalreversible}) holds for
$A=N_i(\mathbf{c})$ and $B=N_j(\mathbf{e})$ which reduces to
%
%e3.11 #&#
\begin{equation}
\label{e22} \mu(i)\mu_i(\mathbf{c})\nu_{i,j}(
\mathbf{e}_i)\frac{c_j }{|\mathbf{c}|} =\mu(j)\mu_j(\mathbf{e})
\nu_{j,i}(\mathbf{c}_j)\frac{e_i }{|\mathbf{e}|}.
\end{equation}
Combining (\ref{e2}) and (\ref{e22}), we get
%
%e3.12 #&#
\begin{equation}
\label{e12} \frac{\nu_{i,j}(\mathbf{d}_i) }{
\nu_{i,j}(\mathbf{e}_i)}=\frac{\mu_j(\mathbf{d}){d_i }/{|\mathbf{d}|} }{
\mu_j(\mathbf{e}){e_i }/{|\mathbf{e}|}}
\end{equation}
holds for all $N_j(\mathbf{d}),N_j(\mathbf{e})$ for which $\mu_j(\mathbf{d}),\mu_j(\mathbf{e})>0$ with $d_i,e_i>0$.

Now for all $k$, $i_1,\dots,i_k=i_1\in[n]$ and $\mathbf{c}(1),\dots,\mathbf{c}(k)=\mathbf{c}(1)$ such that $\mu_{i_p}
(\mathbf{c}(p))>0$ with $c(p)_{i_{p-1}},c(p)_{i_{p+1}}>0$ for
$p=1,\dots,k-1$, we have
\[
\prod_{p=1}^{k-1}\frac{\mu_{i_p}(\mathbf{c}(i))
}{\mu_{i_{p+1}}(\mathbf{c}(i+1))}=1.
\]
Using (\ref{e12}), we have
%
%e3.13 #&#
\begin{equation}
\label{ecycles} \prod_{p=1}^{k-1}
\frac{\nu_{i_p,i_{p+1}} (\mathbf{c}(p+1)_{i_p} )
}{\nu_{i_{p+1},i_p} (\mathbf{c}(p)_{i_{p+1}} )} \frac{c(p)_{i_{p+1}} }{ c(p+1)_{i_p}}=1.
\end{equation}

In Lemma~\ref{lcycles}, we prove that (\ref{ecycles}) holds for all
$i_1,\dots,i_k$
and $\mathbf{c}(1),\dots,\mathbf{c}(k)$ as above if it holds only for a
special family.
Then in Theorem~\ref{t2}, we prove that if the conditional offspring
distributions of
$\nu$ satisfy (\ref{ecycles}) for this special family, as well as
(\ref{eijji}), then there exists a
measure $\mu$ such that $\mu$ is reversible for simple random walk
and $\mu\sim\nu$.

%%%%%%%%%%%%%%%%%%%%
%%% %%%
%%% LEMMA CYCLES %%%
%%% %%%
%%%%%%%%%%%%%%%%%%%%

%le3.4 #&#
\begin{lemma}\label{lcycles}
Let $n \in\N$. Let $\nu$ be the distribution of a
multi-type Galton--Watson tree whose types are ordered pairs of
elements of $[n]$. Then (\ref{ecycles}) holds for all $k$,
$i_1,\dots,i_k=i_1\in[n]$ and $\mathbf{c}(1),\dots,\mathbf{c}(k)=\mathbf{c}(1) \in\mathbf{S}_n$ such that $\nu_{i_p,i_{p+1}}
(\mathbf{c}(p+1)_{i_p} )>0$ for $p=1,\dots,k-1$ if and only if
(\ref{ecycles}) holds under the restriction that $\mathbf{c}(s) \not=
\mathbf{c}(t)$ whenever $i_s=i_t$ for $s,t \in[k-1]$.
\end{lemma}

%%%%%%%%%%%%%%%%%%%%%%%%%%%%%
%%% PROOF OF LEMMA CYCLES %%%
%%%%%%%%%%%%%%%%%%%%%%%%%%%%%

\begin{pf}
The latter condition is a special case of the former so we
prove it implies the former statement. Let $i_1,\dots,i_m=i_1\in
[n]$ and $\mathbf{d}(1),\dots,\mathbf{d}(m)=\mathbf{d}(1)$ such that
$\nu_{i_p,i_{p+1}} (\mathbf{d}(p+1)_{i_p} )>0$ for
$p=1,\dots,m-1$.

If for some $l,u \in[m]$, $\mathbf{d}(l),\dots,\mathbf{d}(u)=\mathbf{d}(l)\in
\mathbf{S}_n$ satisfy the latter condition in the statement of the
lemma, then
\begin{eqnarray*}
&&\prod_{p=1}^{m-1} \frac{\nu_{i_p,i_{p+1}} (\mathbf{d}(p+1)_{i_p} )
}{\nu_{i_{p+1},i_p} (\mathbf{d}(p)_{i_{p+1}} )}
\frac{d(p)_{i_{p+1}} }{ d(p+1)_{i_p}}
\\
&&\quad=\prod_{p=l}^{u-1} \frac{\nu_{i_p,i_{p+1}} (\mathbf{d}(p+1)_{i_p} ) }{
\nu_{i_{p+1},i_p} (\mathbf{d}(p)_{i_{p+1}} )}
\frac{d(p)_{i_{p+1}}
}{ d(p+1)_{i_p}}\prod_{p\in[1,l-1] \cup[u,m-1]} \frac{\nu_{i_p,i_{p+1}} (\mathbf{d}(p+1)_{i_p} ) }{
\nu_{i_{p+1},i_p} (\mathbf{d}(p)_{i_{p+1}} )}
\frac{d(p)_{i_{p+1}}
}{ d(p+1)_{i_p}}
\\
&&\quad=\prod_{p\in[1,l-1] \cup[u,m-1]} \frac{\nu_{i_p,i_{p+1}} (\mathbf{d}(p+1)_{i_p} ) }{
\nu_{i_{p+1},i_p} (\mathbf{d}(p)_{i_{p+1}} )} \frac{d(p)_{i_{p+1}}
}{ d(p+1)_{i_p}}
\end{eqnarray*}
since the first product in the middle line is
equal to 1. Now $i_1,\dots,i_{l-1},i_u,\dots,i_m$ and $\mathbf{d}(1),\dots,\mathbf{d}(l-1),\mathbf{d}(u),\dots,\mathbf{d}(m)=\mathbf{d}(1)$
satisfy the former condition. We can repeat the same process until
$\mathbf{d}(s) \not= \mathbf{d}(t)$ whenever $i_s=i_t$ for $s,t \in
[k-1]$. Since at each step, the product is preserved, the result
follows.
\end{pf}

%%%%%%%%%%%%%%%%%%%%
%%% MAIN THEOREM %%%
%%%%%%%%%%%%%%%%%%%%

%th3.5 #&#
\begin{theorem}\label{t2}
Let $n \in\N$. Let $\nu$ be a strongly connected
multi-type Galton--Watson measure whose types are ordered pairs of
elements of $[n]$. Then there exists a measure $\mu$ on $\mathbf{T}_n$
such that $\mu$ is reversible for simple random walk and $\mu\sim
\nu$ if and only if the following two conditions hold.
\begin{longlist}[(ii)]

\item [(i)] For all $i,j,k\in[n]$, $\mathbf{c} \in\mathbf{S}_n$ such that
$\nu_{j,i}(\mathbf{c}_j)>0$ and $c_k>0$, we have $\nu_{k,i}(\mathbf{c}_k)>0 ,$

\item [(ii)] (\ref{ecycles}) holds for all $k$, $i_1,\dots,i_k=i_1\in[n]$ and
$\mathbf{c}(1),\dots,\mathbf{c}(k)=\mathbf{c}(1) \in\mathbf{S}_n$ such that
$\nu_{i_p,i_{p+1}} (\mathbf{c}(p+1)_{i_p} )>0$ for
$p=1,\dots,k-1$ and for all $s,t \in[k-1]$, if $i_s=i_t$, then\break
$\mathbf{c}(s) \not= \mathbf{c}(t)$.
\end{longlist}
\end{theorem}

%%%%%%%%%%%%%%%%%%%%%%%%%%%%%
%%% PROOF OF MAIN THEOREM %%%
%%%%%%%%%%%%%%%%%%%%%%%%%%%%%

\begin{pf}
By Lemma~\ref{lcycles}, it is enough to prove the theorem when in
(ii), (\ref{ecycles}) holds for all cycles.

($\Leftarrow$) Assume (i) and (ii) hold. Since the $\mu$-conditional
distributions of the descendant trees of the neighbors of the root
are determined by $\nu$, in order to define $\mu$, it's enough to
specify the values $\mu(i)$ and $\mu_i(\mathbf{c})$ for all $i\in[n]$
and $\mathbf{c} \in\mathbf{S}_n$.

For all $i\in[n]$ and $\mathbf{c} \in\mathbf{S}_n$, let $\mu_i(\mathbf{c})=0$ if $\nu_{j,i}(\mathbf{c}_j)=0$ for some (equivalently for all by
assumption (i)) $j\in[n]$ such that $c_j>0$. Let $\mu_i(\mathbf{c})=1$
if for all $j \in[n]$ and $\mathbf{d} \in\mathbf{S}_n$ such that
$\nu_{j,i}(\mathbf{d}_j)>0$, we have $\mathbf{d}=\mathbf{c}$.

For all $i\in[n]$ let $F_i:=\{\mathbf{c} \in\mathbf{S}_n\dvt \exists j,
\nu_{j,i}(\mathbf{c}_j)>0\}$. For each $i\in[n]$ define the probabilities
$\{\mu_i(\mathbf{c})\dvt \mathbf{c} \in F_i\}$ to
be the unique solution of the equations
%
%e3.14 #&#
\begin{equation}
\label{e514} \sum_{\mathbf{c}\in F_i}\mu_i(
\mathbf{c})=1 \qquad\mbox{for } i \in[n]
\end{equation}
and
%
%e3.15 #&#
\begin{equation}
\label{e44} \frac{\mu_i(\mathbf{c}) }{\mu_i(\mathbf{d})}=\prod_{p=1}^{m-1}
\frac{\nu_{i_{p+1},i_p} (\mathbf{c}(p)_{i_{p+1}} )
}{\nu_{i_p,i_{p+1}} (\mathbf{c}(p+1)_{i_p} )}\frac{|\mathbf{c}(p)| }{|\mathbf{c}(p+1)|}\frac{c(p+1)_{i_p}
}{ c(p)_{i_{p+1}}}
\end{equation}
for $i_1=i,\dots,i_m=i\in[n]$ and
$\mathbf{c}(1)=\mathbf{c},\dots,\mathbf{c}(m)=\mathbf{d}$ such that
$\nu_{i_p,i_{p+1}} (\mathbf{c}(p+1)_{i_p} )>0$ for
$p=1,\dots,m-1$.

Define $\{\mu(i)\dvt i \in[n]\}$ to be the unique solution of the equations
$\sum_{i \in[n]} \mu(i)=1$ for $i \in[n]$ and
%
%e3.16 #&#
\begin{equation}
\label{e316} \frac{\mu(i) }{\mu(j)}=\frac{\mu_j(\mathbf{d})\nu_{j,i}(\mathbf{e}_j){d_i }/{|\mathbf{d}|} }{\mu_i(\mathbf{e})\nu_{i,j}
(\mathbf{d}_i){e_j }/{|\mathbf{e}|}}
\end{equation}
for $\mathbf{d} \in F_j,\mathbf{e} \in F_i$ with $d_i,e_j>0$.

In (\ref{e44}) since $\nu_{i_p,i_{p+1}} (\mathbf{c}(p+1)_{i_p}
)>0$ for
$p=1,\dots,m-1$, condition (i) in Theorem~\ref{t2} implies that $\nu_{i_{p+1},i_p}
(\mathbf{c}(p)_{i_{p+1}} )>0$ for $p=2,\dots,m-1$. Since $\mathbf{c} \in F_i$,
condition (i) also implies that $\nu_{i_2,i_1}(\mathbf{c}(1)_{i_2})>0$.

The right-hand side of (\ref{e44}) does not depend on
$i_1=i,\dots,i_m=i\in[n]$ and $\mathbf{c}(1)=\mathbf{c},\dots,\mathbf{c}(m)=\mathbf{d}$: If $l_1=i,\dots,l_t=i\in[n]$ and $\mathbf{e}(1)=\mathbf{c},
\dots,\mathbf{e}(t)=\mathbf{d} \in\mathbf{S}_n$ are such that
$\nu_{l_p,l_{p+1}} (\mathbf{e}(p+1)_{l_p} )>0$ for
$p=1,\dots,t-1$, then by Lemma~\ref{lcycles} and condition (ii) of the
theorem $i_1,\dots,i_m=i=l_t,l_{t-1},\dots,l_1$ and
$\mathbf{c}(1),\dots,\mathbf{c}(m)=\mathbf{d}=\mathbf{e}(t),\dots,\mathbf{e}(1)=\mathbf{c}=\mathbf{c}(1)$ satisfy
\[
\prod_{p=1}^{m-1} \frac{\nu_{i_p,i_{p+1}} (\mathbf{c}(p+1)_{i_p} )
}{\nu_{i_{p+1},i_p} (\mathbf{c}(p)_{i_{p+1}}
)}
\frac{c(p)_{i_{p+1}} }{ c(p+1)_{i_p}} \prod_{p=1}^{t-1}
\frac{\nu_{l_{p+1},l_p} (\mathbf{e}(p)_{l_{p+1}} )
}{
\nu_{l_p,l_{p+1}} (\mathbf{e}(p+1)_{l_p} )}\frac{e(p+1)_{l_p} }{
e(p)_{l_{p+1}}}=1,
\]
which gives
\begin{eqnarray*}
&&\prod_{p=1}^{m-1} \frac{\nu_{i_p,i_{p+1}} (\mathbf{c}(p+1)_{i_p} )
}{\nu_{i_{p+1},i_p} (\mathbf{c}(p)_{i_{p+1}}
)}
\frac{c(p)_{i_{p+1}} }{ c(p+1)_{i_p}}\frac{|\mathbf{c}(p+1)| }{
|\mathbf{c}(p)|}\\
&&\quad{}\times \prod_{p=1}^{t-1}
\frac{\nu_{l_{p+1},l_p} (\mathbf{e}(p)_{l_{p+1}} ) }{\nu_{l_p,l_{p+1}} (\mathbf{e}(p+1)_{l_p} )}\frac{e(p+1)_{l_p} }{ e(p)_{l_{p+1}}} \frac{|\mathbf{e}(p)|
}{|\mathbf{e}(p+1)|}=1.
\end{eqnarray*}
Therefore, we have
\[
\prod_{p=1}^{m-1} \frac{\nu_{i_{p+1},i_p} (\mathbf{c}(p)_{i_{p+1}} )
}{\nu_{i_p,i_{p+1}} (\mathbf{c}(p+1)_{i_p} )}
\frac{c(p+1)_{i_p}
}{ c(p)_{i_{p+1}}}\frac{|\mathbf{c}(p)| }{
|\mathbf{c}(p+1)|} = \prod_{p=1}^{t-1}
\frac{\nu_{l_{p+1},l_p} (\mathbf{e}(p)_{l_{p+1}} ) }{\nu_{l_p,l_{p+1}}
(\mathbf{e}(p+1)_{l_p} )}\frac{e(p+1)_{l_p} }{ e(p)_{l_{p+1}}}\frac{|\mathbf{e}(p)|
}{|\mathbf{e}(p+1)|}.
\]

Now, (\ref{e44}) gives the relative weights of the numbers
$\mu_i(\mathbf{c})$ which add to 1 and therefore $\mu_i(\mathbf{e})$ for
$\mathbf{e} \in F_i$ are uniquely determined by (\ref{e44}) and (\ref{e514}).
If $\mathbf{e} \notin F_i$, then $\mu_i(\mathbf{e})=0$. Since $i$ is
arbitrary, $\mu_i(\mathbf{e})$ is defined for all $i$ and $\mathbf{e}$.

The right-hand side of (\ref{e316}) does not depend on $\mathbf{e}$ and
$\mathbf{d}$. We need to show that
%
%e3.17 #&#
\begin{equation}
\label{e317} \frac{\mu_j(\mathbf{d})\nu_{j,i}(\mathbf{e}_j){d_i }/{|\mathbf{d}|} }{\mu_i(\mathbf{e})\nu_{i,j}
(\mathbf{d}_i){e_j }/{|\mathbf{e}|}} =\frac{\mu_j(\mathbf{f})\nu_{j,i}(\mathbf{c}_j){f_i }/{|\mathbf{f}|} }{
\mu_i(\mathbf{c})\nu_{i,j}(\mathbf{f}_i){c_j }/{|\mathbf{c}|}}
\end{equation}
whenever $\mu_j(\mathbf{d}),\mu_j(\mathbf{f}),\mu_i(\mathbf{c}),\mu_j(\mathbf{e})>0$ with $d_i,f_i,c_j,e_j>0$, in which case the
remaining terms are all positive. In fact it is enough to show (\ref
{e317}) when $\mathbf{d}=\mathbf{f}$ since
repeated application of (\ref{e317}) in that case will prove (\ref{e317})
in the general case. When
$\mathbf{d}=\mathbf{f}$, (\ref{e317}) reduces to
%
%e3.18 #&#
\begin{equation}
\label{e4} \frac{\nu_{j,i}(\mathbf{c}_j) }{
\nu_{j,i}(\mathbf{e}_j)}=\frac{\mu_i(\mathbf{c}){c_j }/{|\mathbf{c}|} }{
\mu_i(\mathbf{e}){e_j }/{|\mathbf{e}|}}
\end{equation}
which is
equivalent to (\ref{e12}) and a special case of (\ref{e44}): If
$\nu_{j,i}(\mathbf{c}),\nu_{j,i}(\mathbf{e})>0$, then let $\mathbf{d} \in
\mathbf{S}_n$ such that $d_i>0$ and $\nu_{i,j}(\mathbf{d}_i)>0$. Then for
$\mathbf{c}(1)=\mathbf{c},\mathbf{c}(2)=\mathbf{d},\mathbf{c}(3)=\mathbf{e}$ and
$i_1=i_3=i,i_2=j$, (\ref{e44}) reduces to (\ref{e4}).

Since $\sum_{i \in[n]}{\mu(i)=1}$, the probabilities $\mu(i)$ for
$i \in[n]$
are uniquely determined. We have that (\ref{e316}) is equivalent to
(\ref{e2}). Since (\ref{e12}) holds, condition (i)
of the theorem implies (\ref{emunu}). By Lemma~\ref{lreversible},
$\mu$ is reversible.

($\Rightarrow$) This direction has been proved by Lemma~\ref{l1},
Lemma~\ref{lcycles} and the discussion before Lemma~\ref{lcycles}.
\end{pf}

%%%%%%%%%%%%%%%%%%%%%%%%%%%%%%%%%%%%%%%%%%
%%%%%%%% EXAMPLE 1 %%%%%%%%%%%%%%%%%%%
%%%%%%%%%%%%%%%%%%%%%%%%%%%%%%%%%%%%%%%%%%

%ex3.6 #&#
\begin{example}\label{xex1}
Let $\mathbf{a}=(1,1),\mathbf{b}=(2,1)$ and $\mathbf{d}=(0,3)$.
Let $\nu$ be a multi-type Galton--Watson measure whose types are
$\{(i,j)\dvt i,j \in\{1,2\}\}$ with offspring distribution as follows
\begin{eqnarray*}
\nu_{1,1}(0,1)&=&\nu_{1,1}(\mathbf{a}_1)=
\tfrac{3 }{7},\qquad \nu_{1,1}(1,1)=\nu_{1,1}(
\mathbf{b}_1)=\tfrac{4 }{
7},
\\
\nu_{2,1}(1,0)&=&\nu_{2,1}(\mathbf{a}_2)=
\tfrac{3 }{5}, \qquad \nu_{2,1}(2,0)=\nu_{2,1}(
\mathbf{b}_2)=\tfrac{2 }{
5},
\\
\nu_{2,2}(0,2)&=&\nu_{2,2}(\mathbf{d}_2)=
\tfrac{2 }{3},\qquad \nu_{2,2}(1,0)=\nu_{2,2}(
\mathbf{a}_2)=\tfrac{1 }{
3},
\\
\nu_{1,2}(0,1)&=&\nu_{1,2}(\mathbf{a}_1)=1.
\end{eqnarray*}

Let's check that $\nu$ satisfies the two conditions of Theorem \ref
{t2}. We
have $\nu_{1,1}(\mathbf{a}_1)$ and $\nu_{2,1}(\mathbf{a}_2)$ are both
positive, exactly as condition (i) asks for: In order to have $\mu
\sim\nu$ and $\nu_{1,1}(\mathbf{a}_1)$ together, one has to have
$\mu(N_1(\mathbf{a}))>0$ which in turn implies $\nu_{2,1}(\mathbf{a}_2)>0$. Similarly $\nu_{1,1}(\mathbf{b}_1)$
and $\nu_{2,1}(\mathbf{b}_2)$ are both positive and so are $\nu_{2,2}(\mathbf{a}_2)$ and
$\nu_{1,2}(\mathbf{a}_1)$. Since $\mathbf{d}_1=0$, condition (i) says
$\nu_{2,2}(\mathbf{d}_2)>0$ implies $\nu_{2,2}(\mathbf{d}_2)>0$, which
certainly holds.

For condition (ii), let us verify (\ref{ecycles}) for a particular
sequence. One actually has to check 4 such sequences in total. Let
$i_1=1,i_2=2,i_3=2,i_4=1,i_5=1$ and $\mathbf{c}(1)=\mathbf{a},\mathbf{c}(2)=\mathbf{a},
\mathbf{c}(3)=\mathbf{a},\mathbf{c}(4)=\mathbf{b},\mathbf{c}(5)=\mathbf{a}$. It might be easier to imagine a cycle of 4 vertices labeled
$i_1=i_5,i_2,i_3,i_4,$ and the vertex labeled $i_s$ has $\mathbf{c}(s)$
neighbors. Starting the cycle at a different vertex, we can see this
is equivalent to showing that (\ref{ecycles}) holds for,
$i_1=1,i_2=1,i_3=2,i_4=2,i_5=1$ and $\mathbf{c}(1)=\mathbf{b},\mathbf{c}(2)=\mathbf{a},\mathbf{c}(3)=\mathbf{a},\mathbf{c}(4)=\mathbf{a},
\mathbf{c}(5)=\mathbf{b}$.

One then can determine
$\mu_1(\mathbf{a}),\mu_1(\mathbf{b}),\mu_2(\mathbf{d}),\mu_2(\mathbf{a})$ by
(\ref{e4}). In fact they are all equal to $\frac{1 }{2}$. Finally, one
can determine $\mu(1)$ and $\mu(2)$ by (\ref{e316}), which are $\frac{3
}{8}$ and $\frac{5 }{8}$ respectively. Since $\mu\sim\nu$, we
have completely determined $\mu$.
\end{example}

Next, we parametrize all reversible measures $\mu$ for simple random
walk such that $\mu\sim\nu$ for some strongly connected
Galton--Watson measure $\nu$. To state our result, Proposition~\ref
{pdimension},
we first prove two lemmas.

%%%%%%%%%%%%%%%%%%%%%%%%%%%%%%%%%%%%%%%%%%%%%%
%%%% LEMMA EXCHANGE PROPERTY OF MU PRIME %%%%
%%%%%%%%%%%%%%%%%%%%%%%%%%%%%%%%%%%%%%%%%%%%%%

%le3.7 #&#
\begin{lemma}\label{lpropertyofmuprime}
Assume that $\mu\sim\nu$ and that $\mu$ is stationary
for simple random walk. Then for all $i,j \in[n]$ we have $i \in A(\mu_j)$ if and only if $j \in A(\mu_i)$.
\end{lemma}

%%%%%%%%%%%%%%%%%%%%%%%%%%%%%%%%%%%%%%%%%%%%%%%%%%%%%%%
%%%% PROOF OF LEMMA EXCHANGE PROPERTY OF MU PRIME %%%%
%%%%%%%%%%%%%%%%%%%%%%%%%%%%%%%%%%%%%%%%%%%%%%%%%%%%%%%
%
\begin{pf}
It is enough to prove one direction since $i$ and $j$ are
interchangeable. Let $i,j \in[n]$ and assume $j \in A(\mu_i)$.
Then there exists $\mathbf{c} \in\mathbf{S}_n$ such
that $c_j>0$ and $\mu_i(\mathbf{c})>0$. Since $\mu(i)>0$, we have $\mu
(N_i(\mathbf{c}))>0$. The argument in the first paragraph of the proof of
Lemma~\ref{l1} gives $\mu(B_{j,i,\mathbf{c}_j})>0$ which implies $i \in
A(\mu_j)$.
\end{pf}

%%%%%%%%%%%%%%%%%%%%%%%%%%%%%%%%%%%%%%%%%%%%
%%%%% LEMMA PROPERTY OF MU PRIME TWO %%%%%%
%%%%%%%%%%%%%%%%%%%%%%%%%%%%%%%%%%%%%%%%%%%%

%le3.8 #&#
\begin{lemma}\label{lpropertyofmuprime2}
Assume that $\mu\sim\nu$ and $\mu$ is stationary
for simple random walk. Then $\nu$ is strongly connected if and only if $\forall i,j \in[n]$, there exist $i=i_1,\dots,i_m=j$ and
$\mathbf{c}(2),\dots, \mathbf{c}(m-1) \in\mathbf{S}_n$ such that
%
%e3.19 #&#
\begin{equation}
\label{emuicondition} c(k)_{i_{k-1}},c(k)_{i_{k+1}}>0 \quad\mbox{and}\quad
\mu_{i_k}\bigl(\mathbf{c}(k)\bigr)>0\qquad \forall k \in\{2,\dots,m-1
\}.
\end{equation}
\end{lemma}

\begin{pf}
Assume $\nu$ is strongly connected. Fix $i,j \in[n]$. Since $\mu
(i),\mu(j)>0$ and $\nu$ is strongly connected,
there exists $i_1=i,\dots,i_m=j$ such that $(i_k,i_{k+1})\in A({\nu
}_{i_{k-1},i_k})$.
Then for each $k=2,\dots,m-1$
there exists $\mathbf{c}(k) \in\mathbf{S}_n$ with
$c_{i_{k+1}},c_{i_{k-1}}>0$ such that $\nu_{i_{k-1},i_k}(\mathbf{c}(k)_{i_{k-1}})>0$.
Lemma~\ref{l1} gives that $\mu_{i_k}(\mathbf{c})>0$.

The other direction can be proved in a similar fashion.
\end{pf}

Let $M$ be a subset of $ \{N_i(\mathbf{c})\dvt i \in[n], \mathbf{c} \in\mathbf{S}_n\}$. If $\mu$ is a reversible measure for simple random walk on
$\mathbf{T}_n$ such that $\mu\sim\nu$ for some strongly connected
Galton--Watson measure $\nu$ and such that the following holds
%
%e3.20 #&#
\begin{equation}
\label{enbr} \forall i \in[n] \mbox{ and } \mathbf{c} \in\mathbf{S}_n
\qquad \mu_i(\mathbf{c})>0 \quad\iff \quad N_i(\mathbf{c}) \in M,
\end{equation}
then by Lemma~\ref{lpropertyofmuprime} and Lemma \ref
{lpropertyofmuprime2} the
elements of $M$ satisfy the following two conditions.
%
%
%e3.21 #&#
\begin{equation}
\begin{tabular}{p{330pt}@{}}
\label{econdition1}  (i) $\quad\forall i,j \in[n]$  we have  $\exists
N_i(\mathbf{c})\in M$  with  $c_j>0 \iff\exists
N_j(\mathbf{d}) \in M$  with  $d_i>0$,
\end{tabular}
\end{equation}
\vspace*{-19pt}
%e3.22 #&#
\begin{equation}
\begin{tabular}{p{330pt}@{}}
\label{econdition2}  (ii) \hspace*{-2.5pt}$\quad \forall i,j \in[n] \, \exists N_{i_1}(
\mathbf{c}(1)),\dots,N_{i_k}(\mathbf{c}(k))\in M$
such that  $i_1=i,i_k=j$  and
$\mathbf{c}(s)_{i_{s-1}}>0$
 and  $\mathbf{c}(s)_{i_{s+1}}>0$  for all  $s \in[k]$.
\end{tabular}
\end{equation}

Let $M_r$ be the set of reversible measures $\mu$ for simple random
walk on $\mathbf{T}_n$ that satisfy (\ref{enbr}) and $\mu\sim\nu$ for
some strongly connected Galton--Watson measure $\nu$.

%%%%%%%%%%%%%%%%%%%%%%%%%%%%%
%%% PROPOSITION DIMENSION %%%
%%%%%%%%%%%%%%%%%%%%%%%%%%%%%

%pr3.9 #&#
\begin{proposition}\label{pdimension}
Reversible measures $\mu$ which satisfy $\mu\sim
\nu$ for some strongly connected Galton--Watson measure $\nu$ can be
parametrized in the following way. Let $M$ be a subset of
$ \{N_i(\mathbf{c})\dvt i \in[n], \mathbf{c} \in\mathbf{S}_n\}$ that satisfies
(\ref{econdition1}) and (\ref{econdition2}) above. Let $\{\mu'_i(\mathbf{c}) \in(0,1)\dvt i \in[n],\mathbf{c} \in\mathbf{S}_n\}$ satisfy (\ref{enbr})
and $ \sum_{\mathbf{c} \in\mathbf{S}_n} \mu'_i(\mathbf{c})=1$ for each $i
\in[n]$. In terms of all the preceding parameters, define $\mu$ by
\begin{longlist}[(iii)]

\item [(iii)] $\mu_i(\mathbf{c})=\mu'_i(\mathbf{c})$ for all $i\in[n],\mathbf{c}\in
\mathbf{S}_n$,

\item [(iv)] conditioned on the labels of the root and its neighbors the
descendant subtrees of the neighbors of the root are independent
multi-type Galton--Watson trees with conditional offspring
distributions $\{\nu_{i,j}(\mathbf{c})\dvt i,j \in[n],\mathbf{c} \in\mathbf{S}_n\}$ that are the unique solutions of the equations
%
%e3.23 #&#
\begin{equation}
\label{e423} \frac{\nu_{i,j}(\mathbf{d}_i) }{
\nu_{i,j}(\mathbf{e}_i)}=\frac{\mu'_j(\mathbf{d}){d_i }/{|\mathbf{d}|} }{
\mu'_j(\mathbf{e}){e_i }/{|\mathbf{e}|}}
\end{equation}
for $N_j(\mathbf{d}),N_j(\mathbf{e})\in M$ with $d_i,e_i>0$ and $ \sum_{N_j(\mathbf{d})
\in M} \nu_{i,j}(\mathbf{d}_i)=1$,

\item [(v)] $\{\mu(i)\dvt i \in[n]\}$ are the unique solutions of the equations
%
%e3.24 #&#
\begin{equation}
\label{e624} \frac{\mu(i) }{\mu(j)}=\frac{\mu'_j(\mathbf{d})\nu_{j,i}(\mathbf{e}_j){d_i }/{|\mathbf{d}|} }{\mu'_i(\mathbf{e})\nu_{i,j}
(\mathbf{d}_i){e_j }/{|\mathbf{e}|}}
\end{equation}
for $N_j(\mathbf{d}),N_i(\mathbf{e})\in M$ with $d_i,e_j>0$.
\end{longlist}
\end{proposition}

\begin{pf}
We first note that since $\mu\sim\nu$, we have that
$\{\mu_i(\mathbf{c})\dvt N_i(\mathbf{c}) \in M\}, \{\mu(i)\dvt i\in[n]\}$ and
$\{\nu_{i,j}(\mathbf{c})\dvt N_i(\mathbf{c}_j) \in M,c_j>0\}$ determine $\mu$
uniquely.

Assume that $\{\mu'_i(\mathbf{c})\dvt i \in[n],\mathbf{c} \in\mathbf{S}_n\}$
are as in Proposition~\ref{pdimension} and $\mu$ is defined by
(iii), (iv) and~(v). Now by (iii) of Proposition~\ref{pdimension} $\mu$ satisfies
(\ref{enbr})
and by (\ref{e624}) we have that $\mu(i),\mu_i(\mathbf{c})$ and
$\nu_{i,j}(\mathbf{c})$ are solutions to (\ref{e2}). By (ii) and
(iii), (\ref{emunu}) holds. Thus, by Lemma~\ref{lreversible}, $\mu$ is
reversible for simple random walk. By Lemma~\ref{lpropertyofmuprime2},
$\nu$ is strongly connected.

On the other hand for $\mu\in M_r$, we already showed that (\ref{e624})
and (\ref{e423}) hold. Therefore, our parametrization is complete.
\end{pf}

We illustrate Proposition~\ref{pdimension} with the following example.

%ex3.10 #&#
\begin{example}\label{xex2}
Let $\mathbf{a},\mathbf{b},\mathbf{d}$ be as in Example~\ref{xex1}. Let
\[
M=\bigl\{N_1(\mathbf{a}),N_1(\mathbf{b}),N_2(
\mathbf{d}),N_2(\mathbf{a})\bigr\}.
\]
It is easy to check that $M$ satisfies (\ref{econdition1}) and (\ref
{econdition2}).
Let $s,t \in(0,1)$ and
\[
\mu_1'(\mathbf{a})=s=1-\mu_1'(
\mathbf{b}) \quad\mbox{and}\quad \mu_2'(\mathbf{a})=t=1-
\mu_2'(\mathbf{d}).
\]
Now, (\ref{enbr}) holds. So by Proposition~\ref{pdimension}, using
(\ref{e423}),
we have
\[
\frac{\nu_{1,1}(0,1) }{\nu_{1,1}(1,1)}=\frac{s }{1-s}\frac{{1 }/{2} }{
{2 }/{3}}.
\]
Since $\nu_{1,1}(0,1)+\nu_{1,1}(1,1)=1$, we get
\[
\nu_{1,1}(0,1)=\frac{3s }{4-s} \quad\mbox{and}\quad \nu_{1,1}(1,1)=
\frac{4-4s
}{4-s}.
\]
Similarly one can calculate
\begin{eqnarray*}
\nu_{2,1}(2,0)&=&\frac{2s }{3-s}\quad \mbox{and} \quad\nu_{2,1}(1,0)=
\frac{3-3s
}{
3-s},
\\
\nu_{2,2}(0,2)&=&\frac{2t }{2-t}\quad \mbox{and}\quad \nu_{2,2}(1,0)=
\frac{2-3t
}{2-t}.
\end{eqnarray*}
Since $M$ contains $N_2(\mathbf{a})$ and $N_2(\mathbf{d})$ but $d_1=0$, by
the last equation in (iv) of Proposition~\ref{pdimension}, we have
\[
\nu_{1,2}(0,1)=\nu_{1,2}(\mathbf{a}_1)=1.
\]
The probabilities $\mu_i(\mathbf{c})$ for $N_i(\mathbf{c})$ are determined
by (iii) and $\mu(1),\mu(2)$ are determined by (\ref{e624}). They turn
out to be
\[
\mu(1)=\frac{3-3s }{6-4s}\quad \mbox{and}\quad \mu(2)=\frac{3 - s }{6-4s}.
\]
When $s=t=\frac{1 }{2}$, then we get the measure of Example~\ref{xex1}.
\end{example}
%
%%%%%%%%%%%%%%%%%%%%%%%%
%%% REMARK DIMENSION %%%
%%%%%%%%%%%%%%%%%%%%%%%%

%re3.11 #&#
\begin{remark}\label{rdimension}
Let $n \in\N$ and $M$ be as in Proposition~\ref{pdimension}. Then
any $\mu\in M_r$ is determined by the
probabilities $\mu_i(\mathbf{c})$. Since $\sum_{N_i(\mathbf{c})\in
M}\mu_i(N_i(\mathbf{c}))=1$ for each $i \in[n]$ we have $M_r$ is an
$(|M|-n)-$dimensional affine set.
\end{remark}

%%%%%%%%%%%%%%%%%%%%%%%%%%%%%%%%%%%%%%%%%%%%%%%
%%%%%%%%%%%%%%%%%%%%%%%%%%%%%%%%%%%%%%%%%%%%%%%
%%%%%%%%%%% %%%%%%%%%%%%
%%%%%%%%%%% %%%%%%%%%%%%
%%%%%%%%%%% SECTION A SPECIAL CASE %%%%%%%%%%%%
%%%%%%%%%%% %%%%%%%%%%%%
%%%%%%%%%%% %%%%%%%%%%%%
%%%%%%%%%%%%%%%%%%%%%%%%%%%%%%%%%%%%%%%%%%%%%%%
%%%%%%%%%%%%%%%%%%%%%%%%%%%%%%%%%%%%%%%%%%%%%%%

%s4 ###
%s4 #&#
\section{A special case}

\setcounter{equation}{24}

In this section, we study the case when there is no relabeling: Let
$g$ be the function which maps $(T,o)\in\mathbf{T}_n$ to the tree
obtained by removing the label of $o$ so that the root has no label
and all the other vertices have the same label as before. We give a
complete description of measures $\mu$ that are reversible and $\mu
\sim_g \nu$ for some strongly connected Galton--Watson measure $\nu$.
We prove the following theorem.

%th4.1 #&#
\begin{theorem}\label{tnorelabeling}
Let $\nu$ be a strongly connected multi-type
Galton--Watson measure whose types are elements of $[n]$. Then there
exists a reversible measure $\mu$ for simple random walk on $\mathbf{T}_n$ such that $\mu\sim_g \nu$ if and only if for all $i\in[n]$
and $d \in D(\nu_i)$, given that the number of children of the root
is $d$, the $ \nu_i$-offspring distribution of the root is
multinomial with parameters $\{p_{i,j}\dvt j \in[n]\}$ that do not
depend on $d$ and that satisfy
%
%e4.25 #&#
\begin{equation}
\label{ecyclesnorelabeling} \prod_{s=1}^{m-1}
\frac{p_{i_s,i_{s+1}} }{ p_{i_{s+1},i_s}}=1
\end{equation}
for all $i_1,\dots,i_m=i_1$ such that $p_{i_s,i_{s+1}}>0$.

In this case, the reversible measure $\mu$ is unique, and $d+1 \in
D(\mu_i)$ if and only if $d \in D(\nu_i)$ for all $i \in[n]$, and
given the root has label $i$, the $\mu_i$-probability that the root
has degree $d+1$ and the $\nu_i$-probability that the root has
degree $d$ are equal. The probabilities of the label of the root are
given by the solution of
%
%e4.26 #&#
\begin{equation}
\label{emuilarinorani} \frac{\mu(i) }{
\mu(j)}=\frac{p_{j,i} }{ p_{i,j}}
\end{equation}
and
$\sum_{k=1}^n \mu(i)=1$.
\end{theorem}

%co4.2 #&#
\begin{corollary}\label{cnorelabeling}
A probability measure $\mu\in\mathbf{T}_n$ is
reversible for simple random walk and $\mu\sim\nu$ for some
multi-type Galton--Watson measure $\nu$ if and only if $\mu$ and
$\nu$ are as described in Theorem~\ref{tnorelabeling}.
\end{corollary}

First, we show how Theorem~\ref{tnorelabeling} can be used to
construct a
class of measures $\mu$ on $\mathbf{T}_3$ such that $\mu$ is reversible
for simple random walk and $\mu\sim_g \nu$ for some multi-type
Galton--Watson measure $\nu$. The process is the analog of Proposition
\ref{pdimension}.

%ex4.3 #&#
\begin{example}\label{xex3}
By Theorem~\ref{tnorelabeling}, if $\mu$ is reversible for
simple random walk and $\mu\sim_g \nu$ for some Galton--Watson
measure $\nu$, then given the label of the root is $i$ and the
degree of the root, each neighbor of the root has label $j$, with
probability, say $p_{i,j}$, independent of the labels of the other
neighbors of the root. Let
$p_{1,1},p_{1,2},p_{1,3},p_{2,1},p_{2,3},p_{3,1},p_{3,2} \in(0,1)$
with \mbox{$p_{1,1}+p_{1,2}+p_{1,3}=1$} and $p_{2,1}+p_{2,3}=1$ and
$p_{3,1}+p_{3,2}=1$. So if the root has label~$1$, then it can have
neighbors labeled 1, 2, 3 and if the root has label 2, then it can
have neighbors labeled~1, 3 and if the root has label 3, then it can
have neighbors labeled 1, 2 with positive probability.

Suppose that (\ref{ecyclesnorelabeling}) holds. The degrees for the
root conditioned on the label of the root and their probabilities
can be chosen arbitrarily. For simplicity, let
$D(\mu_2)=D(\mu_3)=\{5\}$ and $D(\mu_1)=\{2,3\}$, that is, almost
surely when the label or the root is 2 or 3, the root has 5
neighbors, otherwise has 2 or 3 neighbors, say with equal
probability.

By Theorem~\ref{tnorelabeling}, the $\mu$-probability that the root has
label $i$ can be determined by (\ref{emuilarinorani}). Now
$\mu(N_i(\mathbf{c}))$ can be determined for any $i \in\{1,2,3\}$ and
$\mathbf{c} \in\mathbf{S}_n$. The descendant subtrees of the neighbors of
the root, conditioned on their labels, are independent multi-type
Galton--Watson trees, whose offspring distributions are given by Theorem
\ref{tnorelabeling}, namely, a vertex labeled $2$ or $3$ will have 4
children with probability 1. Each child of a vertex labeled $2$ has
label 1 or 3 with probabilities $p_{2,1},p_{2,3}$, respectively,
independent of the labels of the other children. Similarly each
child of a vertex labeled $3$ has label 2 or 3 with probabilities
$p_{3,1},p_{3,2}$ respectively, independent of the labels of the
other children. A vertex labeled $1$ has 1 child or 2 children, with
equal probability and given the number of children, each child has
label 1, 2 or 3 with probabilities $p_{1,1},p_{1,2},p_{1,3}$,
respectively, independent of the labels of the other children.

Let us see what restriction (\ref{ecyclesnorelabeling}) imposes on
the $p_{i,j}$s. Suppose we have chosen $p_{1,1},p_{1,2}$ and
$p_{2,1}$ arbitrarily. Letting $i_1=1,i_2=2,i_3=3,i_4=1$ and using
(\ref{ecyclesnorelabeling}) together with \mbox{$p_{3,2}+p_{3,1}=1$}, one
can solve for $p_{3,2}$ and $p_{3,1}$. Thus the multinomial
parameters for the neighbors of a vertex labeled 3 are determined by
the multinomial parameters for the neighbors of vertices of types 1
and 2.

It is easy to see that this is the only restriction. Just like the
case in the previous section, it will be enough to check (\ref
{ecyclesnorelabeling}) for sequences where no label other than the
first label is repeated. The sequences of length 3 do not impose any
restrictions and the other sequences of length four are actually the
same as the sequence above.
\end{example}

We devote the rest of this section to proving Theorem~\ref{tnorelabeling}.
Assume that $\mu$ is a reversible measure for simple random walk on
$\mathbf{T}_n$ and $\mu\sim_g \nu$. Since $\mu\sim_g \nu$, by
stationarity of $\mu$ we have $d \in D(\nu_i)$ if and only if $d+1
\in D(\mu_i)$. We first prove the following claims in order. For all
$i \in[n]$ and $d \in D(\nu_i)$ given the degree of the root is
$d$, we have that (the offspring distribution of) $\nu_i$ is
multinomial. Given the degree of the root is $d+1$, $\mu_i$ is
multinomial with the same parameters as $\nu_i$ (given the root has
degree $d$). For each $i \in[n]$ and $d \in A(\nu_i)$, the
multinomial parameters for $\nu_i$ given the root has degree $d $
do not depend on $d$. For all $i \in[n]$ and $d \in A(\nu_i)$, the
$\mu_i$-probability that the root has degree $d+1$ and the
$\nu_i$-probability that the root has degree $d$ are equal. Our main
tool is Lemma~\ref{lindependenceoftypes}. Before proving that, we make a
basic observation.

%le4.4 #&#
\begin{lemma}\label{lAmui}
Assume that $\mu$ is reversible and $\mu\sim_g \nu$.
Then for all $i,j \in[n]$ we have
%
%e4.27 #&#
\begin{equation}
\label{e7} i\in A(\nu_j) \quad\iff \quad i \in A(\mu_j)\quad \iff\quad j \in
A(\mu_i) \quad\iff\quad j\in A(\nu_i).
\end{equation}
\end{lemma}

\begin{pf}
The first equivalence (and the last) is the analog of Lemma~\ref{l1}
and can be proven exactly in the same way. For the second
equivalence since $\mu$ is stationary, we have $\mu(B_{i,j})>0$ if
and only if $\mu(B_{j,i})>0$. By definition of $B_{i,j}$ we have $j
\in A(\mu_i)$ if and only if $\mu(B_{i,j})>0$. Thus, (\ref{e7})
holds.
\end{pf}

%%%%%%%%%%%%%%%%%%%%%%%%%%%%%%%%%%%
%%%%%%%%%%%%%%%%%%%%%%%%%%%%%%%%%%%
%%% %%%
%%% LEMMA INDEPENDENCE OF TYPES %%%
%%% %%%
%%%%%%%%%%%%%%%%%%%%%%%%%%%%%%%%%%%
%%%%%%%%%%%%%%%%%%%%%%%%%%%%%%%%%%%

%le4.5 #&#
\begin{lemma}\label{lindependenceoftypes}
Let $\rho$ be a probability measure on
$\mathbf{S}^n$. Let $d \in D(\rho)$. Assume that the following two
conditions hold:
\begin{enumerate}[(ii)]
\item [(i)] $\rho(\mathbf{c})>0 \mbox{ for all } \mathbf{c} \in A(\rho) \mbox{
with } |\mathbf{c}|=d,$

\item  [(ii)] $\forall\mathbf{c} \mbox{ such that } |\mathbf{c}|=d,\rho(\mathbf{c})>0,c_j,c_k>0,$ we have
%
%e4.28 #&#
\begin{equation}
\label{emain} \frac{c_j + \delta_{j,k}}{ c_j+1}\frac{\rho(\mathbf{c}) }{\rho(\mathbf{c}_j^k)} =\frac{c_k+1 }{ c_k + \delta_{j,k}}
\frac{\rho(\mathbf{c}_k^j) }{\rho(\mathbf{c})}.
\end{equation}
\end{enumerate}
Then $\rho(\mathbf{c} \vert |\mathbf{c}|=d)$ is
multinomial.
\end{lemma}

%%%%%%%%%%%%%%%%%%%%%%%%%%%%%%%%%%%%%%%%%%%%
%%% PROOF OF LEMMA INDEPENDENCE OF TYPES %%%
%%%%%%%%%%%%%%%%%%%%%%%%%%%%%%%%%%%%%%%%%%%%

\begin{pf}
We write (\ref{emain}) in a different form, which will be used
later in the proof. Setting $j=s$ and $k=t$ in (\ref{emain}) and
rearranging the terms, we have
%
%e4.29 #&#
\begin{equation}
\label{emainprime} \rho(\mathbf{c}) =\sqrt{\frac{(c_s+1)(c_t+1) }{(c_s +
\delta_{s,t})(c_t + \delta_{s,t})}}\sqrt{\rho
\bigl(\mathbf{c}_t^s\bigr)\rho \bigl(\mathbf{c}_s^t
\bigr)}.
\end{equation}

Let $P$ be a measure on finite lists labeled by elements of $[n]$
with the following property: If $B$ and $C$ are two lists whose
elements have labels $\mathbf{c}$ with $|\mathbf{c}|=d \in D(\rho)$, then
$P(B)=P(C)=\rho(\mathbf{c})/ (|\mathbf{c}|/\prod_{k \in[n]}c_k! )$.
Then $P(B)>0$ if and only if $\mathbf{c} \in A(\rho)$. Fix $d \in
D(\rho)$. Let $P(j\mid\mathbf{c})$ denote the $P$-probability that the
first element of a list of size $d$ is labeled $j$ given that the
other elements have labels $\mathbf{c}$. We shall show that
%
%e4.30 #&#
\begin{equation}
\label{e33} P(j\mid\mathbf{c})=P(j\mid\mathbf{d})
\end{equation}
for all $j,\mathbf{c}, \mathbf{d} \in A(\rho)$ with $|\mathbf{c}|=|\mathbf{d}|=d-1$. It is enough to show that
%
%e4.31 #&#
\begin{equation}
\label{e447} P(j\mid\mathbf{c})=P\bigl(j\mid\mathbf{c}_l^m
\bigr)
\end{equation}
for all $l,m,\mathbf{c} \in A(\rho)$ such that $c_l>0$ and $|\mathbf{c}|=d-1$. In (\ref{e33}), we may assume without loss of generality
that $c_j\geq d_j$. Thus, in (\ref{e447}), we may assume that $m\not=
j$, $m\not= l$. We can also assume that $c_j>0$ since if we want to
show (\ref{e33}) for $\mathbf{c}$ with $c_j=0$ (consequently $d_j=0$),
we can pick $\mathbf{e}$ with $e_j>0$ and show that (\ref{e33}) holds
for $\mathbf{e}$ and $\mathbf{d}$ and then for $\mathbf{e}$ and $\mathbf{c}$,
which together imply that (\ref{e33}) holds for $\mathbf{c}$ and $\mathbf{d}$.

All possible assignments of the labels $\mathbf{c}^j$ to a list of
$|\mathbf{c}|+1$ elements have the same probability, which is
$\rho(\mathbf{c}^j)/ ((|\mathbf{c}|+1)! / \prod_{k \in
[n]}(c_k+\delta_{k,j})! )$. Hence, the $P$-probability that a list
of $|\mathbf{c}|+1$ elements are labeled by $\mathbf{c}^j$ and the first
element of the list has label $j$ is given by
\[
\rho\bigl(\mathbf{c}^j\bigr) \biggl(\frac{|\mathbf{c}|! }{\prod_{k \in[n]}c_k!} \bigg/
\frac{(|\mathbf{c}|+1)! }{\prod_{k \in[n]}(c_k+\delta_{k,j})!} \biggr) =\rho\bigl(\mathbf{c}^j\bigr) \biggl(
\frac{c_j+1 }{|\mathbf{c}|+1} \biggr).
\]
Using this
equation, we have
\[
P(j\mid\mathbf{c})=\frac{({c_j+1 })/({|\mathbf{c}|+1})\rho(\mathbf{c}^j) }{\sum_{k \in A}
({c_k+1 })/({|\mathbf{c}|+1})\rho(\mathbf{c}^k)} =\frac{(c_j+1) \rho(\mathbf{c}^j)
}{\sum_{k \in A}(c_k+1) \rho(\mathbf{c}^k)},
\]
where $A:=A({\rho})$. Similarly
\[
P\bigl(j\mid\mathbf{c}_l^m\bigr)=\frac{(c_j+1-\delta_{j,l})
\rho(\mathbf{c}_l^{m,j})}{\sum_{k \in
A}(c_k+1+\delta_{k,m}-\delta_{k,l}) \rho(\mathbf{c}_l^{m,k})}.
\]
Hence,
\[
P(j\mid\mathbf{c})=P\bigl(j\mid\mathbf{c}_l^m\bigr)\quad \iff\quad
\frac{(c_j+1) \rho(\mathbf{c}^j)
}{
\sum_{k \in A}(c_k+1) \rho(\mathbf{c}^k)}= \frac{(c_j+1-\delta_{j,l})
\rho(\mathbf{c}_l^{m,j})}{\sum_{k \in
A}(c_k+1+\delta_{k,m}-\delta_{k,l}) \rho(\mathbf{c}_l^{m,k})},
\]
which is equivalent to
%
%e4.32 #&#
\begin{eqnarray}
\label{e5} &&(c_j+1) \rho\bigl(\mathbf{c}^j\bigr)\sum
_{k \in A}(c_k+1+\delta_{k,m}-
\delta_{k,l}) \rho\bigl(\mathbf{c}_l^{m,k}
\bigr)\nonumber\\[-8pt]\\[-8pt]
&&\quad=(c_j+1-\delta_{j,l}) \rho\bigl(\mathbf{c}_l^{m,j}
\bigr)\sum_{k \in
A}(c_k+1) \rho\bigl(
\mathbf{c}^k\bigr).\nonumber
\end{eqnarray}

We shall show that $\forall k \in A$
%
%e4.33 #&#
\begin{equation}
\label{e6} (c_j+1) \rho\bigl(\mathbf{c}^j\bigr)
(c_k+1+\delta_{k,m}-\delta_{k,l}) \rho\bigl(
\mathbf{c}_l^{m,k}\bigr)=(c_j+1-
\delta_{j,l}) \rho\bigl(\mathbf{c}_l^{m,j}\bigr)
(c_k+1)\rho\bigl(\mathbf{c}^k\bigr).
\end{equation}
This implies
(\ref{e5}), which completes the proof.

Since $c_l>0$ and $c_j>0$, by (\ref{emainprime}), we have the
following equations
\begin{eqnarray*}
\rho\bigl(\mathbf{c}^j\bigr)&=&\sqrt\frac{(c_j+2)(c_l+1+\delta_{l,j}) }{
(c_j+1+\delta_{l,j})(c_l+2\delta_{l,j})}\sqrt
{\rho\bigl(\mathbf{c}_l^{j,j}\bigr)\rho\bigl(
\mathbf{c}^l\bigr)}
\\
&=&\sqrt\frac{(c_j+2-\delta_{l,j})(c_l+1)}{(c_j+1)(c_l+\delta_{l,j})}\sqrt{\rho\bigl(\mathbf{c}_l^{j,j}
\bigr)\rho\bigl(\mathbf{c}^l\bigr)},
\\
\rho\bigl(\mathbf{c}_l^{m,k}\bigr)&=&\sqrt
\frac{(c_m+2+\delta_{m,k})(c_k+2+\delta_{m,k}-\delta_{l,k})}{(c_m+1+2\delta_{m,k})
(c_k+1+2\delta_{m,k}-\delta_{l,k})}\sqrt{\rho\bigl(\mathbf{c}_l^{m,m}
\bigr)\rho \bigl(\mathbf{c}_l^{k,k}\bigr)}
\\
&=&\sqrt\frac{(c_m+2)(c_k+2-\delta_{l,k})}{(c_m+1+\delta_{m,k})(c_k+1+\delta_{m,k}-\delta_{l,k})} \sqrt{\rho\bigl(\mathbf{c}_l^{m,m}
\bigr)\rho\bigl(\mathbf{c}_l^{k,k}\bigr)},
\\
\rho\bigl(\mathbf{c}^k\bigr)&=&\sqrt\frac{(c_k+2)(c_l+1+\delta_{l,k}) }{
(c_k+1+\delta_{l,k})
(c_l+2\delta_{l,k})}\sqrt
{\rho\bigl(\mathbf{c}_l^{k,k}\bigr)\rho\bigl(
\mathbf{c}^l\bigr)}
\\
&=&\sqrt\frac{(c_k+2-\delta_{l,k})(c_l+1) }{(c_k+1)(c_l+\delta_{l,k})} \sqrt{\rho\bigl(\mathbf{c}_l^{k,k}
\bigr)\rho\bigl(\mathbf{c}^l\bigr)},
\\
\rho\bigl(\mathbf{c}_l^{m,j}\bigr)&=&\sqrt
\frac{(c_m+2+\delta_{m,j})(c_j+2+\delta_{m,j}-\delta_{l,j})}{
(c_m+1+2\delta_{m,j})(c_j+1+2\delta_{m,j}-\delta_{l,j})}\sqrt{\rho \bigl(\mathbf{c}_l^{m,m}
\bigr)\rho\bigl(\mathbf{c}_l^{j,j}\bigr)}
\\
&=&\sqrt\frac{(c_m+2)(c_j+2-\delta_{l,j})}{(c_m+1)(c_j+1-\delta_{l,j})} \sqrt{\rho\bigl(\mathbf{c}_l^{m,m}
\bigr)\rho\bigl(\mathbf{c}_l^{j,j}\bigr)}.
\end{eqnarray*}
Substituting these into (\ref{e6}), after cancelations and squaring
the equation, we get
\[
(c_j+1) (c_k+1+\delta_{k,m}-
\delta_{k,l}) (c_l+\delta_{l,k})
(c_m+1) =(c_k+1) (c_j+1-\delta_{j,l})
(c_l+\delta_{l,j}) (c_m+1+\delta_{m,k}).
\]

Suppose $k=j$, then $\delta_{k,m}=\delta_{j,m}=0$. Substituting $j$
for $k$ and using that $\delta_{j,m}=0$, we see that the equation
above holds when $j=k$ and so does (\ref{e6}).

So assume $k \not= j$.

If $k=l$, then $l \not= j$ and $m \not= k$ and we can evaluate all
the $\delta$'s to see that the equation above holds.

If $k\not=l$, the equation above reduces to
\[
(c_j+1) (c_l)=(c_j+1-\delta_{j,l})
(c_l+\delta_{l,j}),
\]
which holds whether $j=l$ or $j \not= l$.
\end{pf}

%%%%% QED %%%%%%%%%%

We shall apply Lemma~\ref{lindependenceoftypes} with $\rho:=\nu_i$
conditioned on the number of children (for each $i$). We now show
that for each $i$, $\nu_i$ satisfies condition (i) of Lemma \ref
{lindependenceoftypes}.

%%%%%%%%%%%%%%%%%%%%%%%%%%%%%%%%%%%
%%% LEMMA TEMELDEGISIMOZELLIGI %%%
%%%%%%%%%%%%%%%%%%%%%%%%%%%%%%%%%%%

%le4.6 #&#
\begin{lemma}\label{ltemeldegisimozelligi}
Assume $\mu$ is stationary for simple
random walk and $\mu\sim_g \nu$.
Suppose $i,l \in[n]$ and $\mathbf{d} \in\mathbf{S}_n$ are such that
$\nu_i(\mathbf{d})>0$ with $d_l>0$. Then $\nu_i(\mathbf{d}_l^m)>0$ for
every $m \in A(\nu_i)$.
\end{lemma}

\begin{pf}
Fix $\mathbf{d},i,l,m$ as in the lemma. Since $m
\in A(\nu_i)$, (\ref{e7}) gives $\mu(B_{m,i})>0$. Since $\mu
\sim\nu$, $\mu(B_{m,i})>0$ and $\nu_i(\mathbf{d})>0$ together imply
$\mu(B_{m,i,\mathbf{d}})>0$. For every $(T,o)\in B_{m,i,\mathbf{d}}$, a
random walker starting at $(T,o)$ moves to a tree $(T,o') \in
B_{l,i,\mathbf{d}_l^m}$ in two steps with positive probability. Since
$\mu$ is stationary and $\mu(B_{m,i,\mathbf{d}})>0$, we have
$\mu(B_{l,i,\mathbf{d}_l^m})>0$. Since $\mu\sim_g \nu$, this implies
$\nu_i(\mathbf{d}_l^m)>0$.
\end{pf}

Now if $\mu$ is reversible for simple random walk and $\mu\sim_g
\nu$ we have that for all $i\in[n]$ condition (i) of Lemma \ref
{lindependenceoftypes}
for $\rho=\nu_i$ is a corollary of Lemma \ref
{ltemeldegisimozelligi}. Next, we show that under the same hypothesis
condition (ii) of Lemma~\ref{lindependenceoftypes} holds for $\rho=\nu_i$
for all $i\in[n]$.

%%%%%%%%%%%%%%%%%%%%%%%%%%%%%%%%%%%%%%%%%%%%%%%%%%%%%%%
%%% NU_I'IN IKINCI KOSULU SAGLADIGININ GOSTERILMESI %%%
%%%%%%%%%%%%%%%%%%%%%%%%%%%%%%%%%%%%%%%%%%%%%%%%%%%%%%%

Assume that $\nu_i(\mathbf{c})>0$ with $c_j,c_k>0$. Then $j,k \in
A(\nu_i)$, which implies $i \in A(\nu_j) \cap A(\nu_k)$ by Lemma
\ref{lAmui}. Let $\mathbf{d},\mathbf{e}$ be such that $d_i,\nu_i(\mathbf{d})>0$
and $\mu_i(\mathbf{e}),e_i>0$. Let $i_1=i,i_2=j,i_3=i,i_4=k,i_5=i$ and
$\mathbf{c}(1)=\mathbf{c}(5)=\mathbf{c}^k,\mathbf{c}(2)=\mathbf{d},\mathbf{c}(3)=\mathbf{c}^j,
\mathbf{c}(4)=\mathbf{e}$. Then since $\mu$ is reversible, (\ref{ecycles})
holds for $i_1,\dots,i_5$ and $\mathbf{c}(1),\dots,\mathbf{c}(5)$ where $\nu_{s,t}=\nu_t$ for all $s,t$
(since~\mbox{$f=g$}),
\[
\frac{\nu_j(\mathbf{d}_i) }{\nu_i(\mathbf{c}^k_j)}\frac{c_j+\delta_{j,k} }{
d_i}\frac{\nu_i(\mathbf{c}) }{\nu_j(\mathbf{d}_i)}\frac{d_i }{ c_j+1}
\frac{\nu_k(\mathbf{e}_i) }{\nu_i(\mathbf{c}^j_k)}\frac{c_k+\delta_{j,k} }{
e_i} \frac{\nu_i(\mathbf{c}) }{\nu_k(\mathbf{e}_i)}\frac{e_i }{
c_k+1}=1,
\]
which is equivalent to (\ref{emain}).

Therefore, we can apply Lemma~\ref{lindependenceoftypes} to $\nu_i$
for $i
\in[n]$ to deduce that conditioned on the degree of the root and
the type of the root is $i$,
the $\nu_i$ offspring distribution of the root is multinomial.

%%%%%%%%%%%%%%%%%%%%%%%%%%%%%%%%%%%%%%%%%%%%%%%%%%%%%%%%%%%%%%
%%% $mu_i$ IN INDEPENDENCEOFTYPES'IN SARTLARINI SAGLAMASI %%%
%%%%%%%%%%%%%%%%%%%%%%%%%%%%%%%%%%%%%%%%%%%%%%%%%%%%%%%%%%%%%%

We next show that Lemma~\ref{lindependenceoftypes} can be applied to
$\mu_i$ for each $i\in[n]$.

%%%%%%%%%%%%%%%
%%% LEMMA 8 %%%
%%%%%%%%%%%%%%%

%le4.7 #&#
\begin{lemma}\label{l8}
Assume $\mu$ is reversible for simple random walk and
$\mu\sim_g \nu$. Suppose $i,l \in[n]$ and $\mathbf{c} \in\mathbf{S}_n$
are such that $\mu_i(\mathbf{c})>0$ with $c_l>0$. Then $\mu_i(\mathbf{c}_l^m)>0$ for every $m \in A(\mu_i).$
\end{lemma}

\begin{pf}
Assume $|\mathbf{c}|=1$. Then $\mu(B_{\mathbf{},i,l})>0$, which
implies $\mu(B_{l,i,\mathbf{}})>0$. Hence $\nu_i(\mathbf{0})>0$. Since
$m\in A(\mu_i)$, $i\in A(\mu_m)$. Consequently the set of trees
whose root has label $m$ and has a neighbor labeled $i$ that does
not have any other neighbor has positive $\mu$-probability. By
stationarity of $\mu$, $\mu(B_{\mathbf{},i,m})>0$. Therefore, we have
$\mu_i(\mathbf{c}_l^m)>0$.

If $|\mathbf{c}|>1$, let $\mathbf{d}:=\mathbf{c}_j$ for some $j$ such that $d_l>0$.
We have $\mu(B_{\mathbf{d},i,j})=\mu (N_i(\mathbf{c}) )>0$. Then
stationarity
of $\mu$ implies that $\mu(B_{j,i,\mathbf{d}})>0$. Since $\mu\sim\nu$,
we have $\nu_i(\mathbf{d})>0$. By Lemma~\ref{ltemeldegisimozelligi},
$\nu_i(\mathbf{d}^m_l)>0$ for all $m\in
A(\nu_i)$. We also have $A(\nu_i)=A(\mu_i)$. Hence for all
$m\in A(\mu_i)$, we have $\nu_i(\mathbf{d}^m_l)>0$, which implies
$\mu(B_{j,i,\mathbf{d}^m_l})>0$. Since $\mu$ is stationary, we get $\mu
(N_i({\mathbf{c}^m_l}) )=\mu(B_{\mathbf{d}^m_l,i,j})>0$.
\end{pf}

As a corollary to Lemma~\ref{l8}, if $\mu$ is reversible for simple
random walk and $\mu\sim_g \nu$, then $\mu_i$ satisfies condition
(i) of Lemma~\ref{lindependenceoftypes}. We now show that in this case
$\mu_i$ also satisfies condition (ii) in the statement
of Lemma~\ref{lindependenceoftypes}.

Since the $\nu_i$ conditional offspring distribution of the root is
multinomial, for all $i \in[n]$ and for all $d \in D(\nu_i)$ there
exist probabilities $\{p_{i,j}^d\dvt j \in[n]\}$ such that for all
$\mathbf{c} \in A(\nu_i)$ with $|\mathbf{c}|=d$, we have
$\nu_i(\mathbf{c})=\frac{|\mathbf{c}|! }{\prod_{j \in[n]} c_j!}\prod_{j\in
[n]}(p^d_{i,j})^{c_j}.$

Fix $\mathbf{c} \in\mathbf{S}_n$ and $i,j,k \in[n]$ such that
$\mu_i(\mathbf{c})>0$ with $c_j,c_k>0$ and $j \not= k$. Since $j\in
A(\mu_i)$, we have $i\in A(\mu_j)$ by Lemma~\ref{lAmui}. Let $\mathbf{d}$ be
such that $\mu_j(\mathbf{d})>0$ and $d_i>0$. Then since $\mu$ is
reversible, (\ref{e2}) holds and is equivalent to
%
%e4.34 #&#
\begin{equation}
\label{e4gg} \mu(i)\mu_i(\mathbf{c})\nu_j(
\mathbf{d}_i)\frac{c_j }{|\mathbf{c}|}=\mu(j)\mu_j(\mathbf{d})
\nu_i(\mathbf{c}_j)\frac{d_i }{|\mathbf{d}|}.
\end{equation}
Similarly we have
\[
\mu(i)\mu_i\bigl(\mathbf{c}^j_k\bigr)
\nu_j(\mathbf{d}_i)\frac{c_j+1 }{|\mathbf{c}^j_k|}=\mu(j)
\mu_j(\mathbf{d})\nu_i(\mathbf{c}_k)
\frac{d_i }{|\mathbf{d}|}.
\]
Combining the two equations, we get
%
%e4.35 #&#
\begin{equation}
\frac{\mu_i(\mathbf{c}) }{\mu_i(\mathbf{c}^j_k)}=\frac{c_j+1
}{ c_j}\frac{\nu_i(\mathbf{c}_j) }{\nu_i(\mathbf{c}_k)}.
\end{equation}
Since
$\nu_i$ is multinomial given the number of children and $d:=|\mathbf{c}_j|=|\mathbf{c}_k|$, we have
%
%e4.36 #&#
\begin{equation}
\label{e16} \frac{\mu_i(\mathbf{c}) }{\mu_i(\mathbf{c}^j_k)}=\frac{c_j+1 }{ c_j}\frac{\nu_i(\mathbf{c}_j) }{\nu_i(\mathbf{c}_k)}=
\frac{c_j+1 }{ c_j}\frac{c_j }{ c_k}\frac{p^d_{i,k} }{
p^d_{i,j}}=\frac{c_j+1 }{ c_k}
\frac{p^d_{i,k} }{ p^d_{i,j}}.
\end{equation}
Interchanging $j$ and $k$ in (\ref{e16}), we get
%
%e4.37 #&#
\begin{equation}
\label{e17} \frac{\mu_i(\mathbf{c}) }{\mu_i(\mathbf{c}^k_j)}=\frac{c_k+1 }{ c_j}\frac{p^d_{i,j}
}{
p^d_{i,k}}.
\end{equation}
Combining (\ref{e16}) and (\ref{e17}), we obtain
\[
\frac{c_k }{ c_{j+1}}\frac{\mu_i(\mathbf{c}) }{\mu_i(\mathbf{c}^j_k)}=\frac{p^d_{i,k} }{
p^d_{i,j}}=\frac{c_k+1 }{ c_j}
\frac{\mu_i(\mathbf{c}^k_j) }{\mu_i(\mathbf{c})},
\]
which is equivalent to (\ref{emain}) when $k \not= j$. When
$j=k$, (\ref{emain}) holds trivially.

Thus for all $i\in[n]$ and for all $d+1 \in D(\mu_i)$, given the
degree of the root is $d+1$ and the type of the root is $i$, the
$\mu_i$-offspring distribution of the root is multinomial.
By (\ref{e17}) the multinomial parameters for $\mu_i$ given that the degree
of the root is $d+1$ and for $\nu_i$ given that the degree of the
root is $d$ are identical.

%le4.8 #&#
\begin{lemma}\label{l9}
Assume $\mu$ is reversible for simple random walk and
$\mu\sim_g \nu$. Then $p_{i,j}^d=p_{i,j}^{d'}$ for all $i \in[n], j
\in A(\nu_i)$ and $d,d' \in D(\nu_i)$.
\end{lemma}

\begin{pf}
Let $i,j,d,d'$ be as in the lemma and $\mathbf{c},\mathbf{e},\mathbf{d}$ be such that $\mu_i(\mathbf{c}),\mu_i(\mathbf{e}),\mu_j(\mathbf{d})>0$
with $c_j,e_j,d_i>0$ and $|\mathbf{c}|=d+1,|\mathbf{e}|=d'+1$. Writing
(\ref{e4gg}) for $\mathbf{e}$ and $\mathbf{d}$, we get
%
%e4.38 #&#
\begin{equation}
\label{e4g} \mu(i)\mu_i(\mathbf{e})\nu_j(
\mathbf{d}_i)\frac{e_j }{|\mathbf{e}|}=\mu(j)\mu_j(\mathbf{d})
\nu_i(\mathbf{e}_j)\frac{d_i }{|\mathbf{d}|}.
\end{equation}
Combining (\ref{e4gg}) and (\ref{e4g}), we get
%
%e4.39 #&#
\begin{equation}
\label{e10} \frac{\mu_i(\mathbf{c})c_j }{\nu_i(\mathbf{c}_j)|\mathbf{c}|}=\frac{\mu_i(\mathbf{e})e_j }{\nu_i(\mathbf{e}_j)|\mathbf{e}|}.
\end{equation}
Since the $\mu$-probability that the root has degree $d+1$ given
that it has label $i$ is equal to $\sum_{|\mathbf{f}|=d+1}\mu_i(\mathbf{f}) ,$ we can rewrite the left-hand side of (\ref{e10}) as
\[
\frac{\sum_{|\mathbf{f}|=d+1}\mu_i(\mathbf{f})\prod_{k \in[n]}
(p_{i,k}^d)^{c_k}{|\mathbf{c}|! }/({\prod_k c_k!})c_j
}{\sum_{|\mathbf{f}|=d}\nu_i(\mathbf{f}) \prod_{k \in[n]}
(p_{i,k}^d)^{c_k-\delta_{k,j}}{(|\mathbf{c}|-1)! }/({\prod_k
(c_k-\delta_{k,j})!})|\mathbf{c}|}= \frac{\sum_{|\mathbf{f}|=d+1}\mu_i(\mathbf{f})
}{\sum_{|\mathbf{f}|=d}\nu_i(\mathbf{f})}p_{i,j}^d.
\]
The right-hand
side of (\ref{e10}) can be calculated in a similar way. Then (\ref
{e10}) reduces to
%
%e4.40 #&#
\begin{equation}
\label{e11} \frac{\sum_{|\mathbf{f}|=d+1}\mu_i(\mathbf{f}) }{\sum_{|\mathbf{f}|=d}\nu_i(\mathbf{f})}p_{i,j}^d=
\frac{\sum_{|\mathbf{f}|=d'+1}\mu_i(\mathbf{f}) }{\sum_{|\mathbf{f}|=d'}\nu_i(\mathbf{f})}p_{i,j}^{d'}
\end{equation}
and (\ref{e11})
holds for all $j \in A(\mu_i)$. Since $\sum_{j \in
A(\mu_i)}p_{i,j}^d=\sum_{j \in A(\mu_i)}p_{i,j}^{d'}=1$, (\ref{e11})
implies that $p_{i,j}^d=p_{i,j}^{d'}:=p_{i,j}$ for all $j \in
A(\mu_i)$.
\end{pf}

We now finish proving the claims in the first
paragraph. Since
\[
\sum_{d+1 \in D(\mu_i)}\sum_{|\mathbf{f}|=d+1}
\mu_i(\mathbf{f})=\sum_{d \in D(\nu_i)}\sum
_{|\mathbf{f}|=d}\nu_i(\mathbf{f})=1,
\]
(\ref{e11}) implies that $\sum_{|\mathbf{f}|=d+1}\mu_i(\mathbf{f})=\sum_{|\mathbf{f}|=d}\nu_i(\mathbf{f})=:p_i^d$ for
all $d \in A(\nu_i)$.

Next, we show (\ref{ecyclesnorelabeling}). The left-hand side of (\ref
{e4gg}) is equal to
\[
\mu(i)p_i^d\prod_{k \in[n]}
\bigl(p_{i,k}^d\bigr)^{c_k}\frac{|\mathbf{c}|! }{\prod_{k \in[n]} c_k!}p_j^{d'}
\prod_{k \in[n]} \bigl(p_{j,k}^{d'}
\bigr)^{d_k-\delta_{i,k}}\frac{d'! }{\prod_k (e_k-\delta_{i,k})!}\frac{c_j }{|\mathbf{c}|}.
\]
Similarly the right-hand side is equal to
\[
\mu(j)p_j^{d'}\prod_{k \in[n]}
\bigl(p_{j,k}^{d'}\bigr)^{d_k}\frac{|\mathbf{d}|! }{
\prod_{k \in[n]} d_k!}p_i^d
\prod_{k \in[n]} \bigl(p_{i,k}^d
\bigr)^{c_k-\delta_{k,j}}\frac{d! }{\prod_k
(c_k-\delta_{k,j})!}\frac{d_i }{|\mathbf{d}|}.
\]
Setting the two equal
to each other we get (\ref{emuilarinorani}). Since for any
$i_1,\dots,i_m=i_1$ we have $\prod_{s=1}^{m-1} \frac{\mu(i_s) }{
\mu(i_{s+1})}=1$, (\ref{emuilarinorani}) gives (\ref{ecyclesnorelabeling}).

Conversely assume that given the degree of the root the
$\nu_i$-offspring distributions of the root are multinomial with
parameters $\{p_{i,j}\dvt j \in[n]\}$ that do not depend on $d$ and
that (\ref{ecyclesnorelabeling}) holds. For all $i \in[n]$ and $d
\in D(\nu_i)$, let $p_i^d$ be the $\nu_i$-probability that the root
has degree $d$.

Define $\mu$ in the following way. Let $\{\mu(i)\dvt  i\in[n]\}$ be the
unique solution of the equations (\ref{emuilarinorani}) and
$\sum_{i=1}^n \mu(i)=1$. For all $i \in[n]$ let the
$\mu_i$-probability that the root has degree $d+1$ given the root
has label $i$ be $p_i^d$. Given the labels of the root and its
neighbors the descendant subtrees of the neighbors of the root are
independent multi-type Galton--Watson trees with conditional
offspring distributions same as conditional offspring distributions
of $\nu$.

Since (\ref{ecyclesnorelabeling}) holds, there is a solution to (\ref
{emuilarinorani}) such that $\mu$ is a probability measure on~$[n]$.
The equations (\ref{emuilarinorani}) give the relative weights of
$\mu(i)$ hence the solution is unique. It is easy to see that (\ref
{e4g}) holds and therefore $\mu$ is reversible.

This completes the proof of Theorem~\ref{tnorelabeling}. We finish
with two
open questions.

What can be said for more general relabeling functions? In other
words are there other Galton--Watson measures $\nu$ for which there
exist a reversible measure $\mu$ for simple random walk and a
relabeling function $f$ (other than the ones we have studied so far)
such that $\mu\sim_f \nu$?

Suppose that when conditioned on the labels of the root and its
neighbors, the descendant subtrees of the root are multi-type
Galton--Watson trees (not necessarily independent). In this case, if
the labels of the root and the neighbors of the root are given, then
for any neighbor of the root, the descendant subtrees of the
neighbors of that neighbor of the root are independent multi-type
Galton--Watson trees. We ask whether reversibility of $\mu$ implies
the conditional independence of the descendant subtrees of the
neighbors of the root.

\section*{Acknowledgement}
I am indebted to Russ Lyons for useful
discussions, suggestions and comments.

\printhistory


\begin{thebibliography}{5}
% BibTex style file: bej.bst, 2011-10-13
% Default style options (sort=1,type=number).
% Used options (sort=1,type=number).

%b1 #&#
\bibitem{aldous}
\begin{barticle}[mr]
\bauthor{\bsnm{Aldous},~\bfnm{David}\binits{D.}} \AND
\bauthor{\bsnm{Lyons},~\bfnm{Russell}\binits{R.}}
(\byear{2007}).
\btitle{Processes on unimodular random networks}.
\bjournal{Electron. J.~Probab.}
\bvolume{12}
\bpages{1454--1508}.
\bid{doi={10.1214/EJP.v12-463}, issn={1083-6489}, mr={2354165}}
\bptok{imsref}%
\end{barticle}
\endbibitem

%b2 #&#
\bibitem{altok}
\begin{barticle}[mr]
\bauthor{\bsnm{Altok},~\bfnm{Serdar}\binits{S.}}
(\byear{2010}).
\btitle{Reversibility of a simple random walk on periodic trees}.
\bjournal{Proc. Amer. Math. Soc.}
\bvolume{138}
\bpages{1101--1111}.
\bid{doi={10.1090/S0002-9939-09-09844-X}, issn={0002-9939}, mr={2566575}}
\bptok{imsref}%
\end{barticle}
\endbibitem

%b3 #&#
\bibitem{augmented}
\begin{barticle}[mr]
\bauthor{\bsnm{Lyons},~\bfnm{Russell}\binits{R.}},
\bauthor{\bsnm{Pemantle},~\bfnm{Robin}\binits{R.}} \AND
\bauthor{\bsnm{Peres},~\bfnm{Yuval}\binits{Y.}}
(\byear{1995}).
\btitle{Ergodic theory on {G}alton--{W}atson trees: Speed of random walk and
dimension of harmonic measure}.
\bjournal{Ergodic Theory Dynam. Systems}
\bvolume{15}
\bpages{593--619}.
\bid{doi={10.1017/S0143385700008543}, issn={0143-3857}, mr={1336708}}
\bptok{imsref}%
\end{barticle}
\endbibitem

%b4 #&#
\bibitem{Takacs}
\begin{barticle}[mr]
\bauthor{\bsnm{Takacs},~\bfnm{Christiane}\binits{C.}}
(\byear{1997}).
\btitle{Random walk on periodic trees}.
\bjournal{Electron. J. Probab.}
\bvolume{2}
\bpages{1--16 (electronic)}.
\bid{doi={10.1214/EJP.v2-15}, issn={1083-6489}, mr={1436761}}
\bptok{imsref}%
\end{barticle}
\endbibitem

\end{thebibliography}
\end{document}